\documentclass{article}
\usepackage{graphicx} 
\usepackage[utf8]{inputenc}
\usepackage[T1]{fontenc}
\usepackage{mathptmx}
\usepackage{etoolbox}
\usepackage{natbib}
\usepackage{relsize}
\usepackage{amssymb}
\usepackage{amsmath}
\usepackage{amsthm}
\usepackage{subcaption}
\usepackage{float}
\usepackage{hyperref}
\usepackage{adjustbox}
\usepackage{xcolor}
\usepackage{soul}
\usepackage{caption}
\usepackage[a4paper, margin=1.5in]{geometry}
\usepackage{graphicx}
\usepackage[percent]{overpic}

\title{Stability of large dynamical systems using Gershgorin Disc Theorem}
\author{
Yuvrajsingh Pravinsingh Patil$^{1}$,
Chandrakala Meena$^{1,*}$
\\[0.5em]
$^{1}$Indian Institute of Science Education and Research (IISER) Pune, India \\
\\[0.5em]
}
\date{}
\begin{document}

\maketitle

\begin{abstract}
Understanding the stability of large nonlinear dynamical systems is a fundamental problem in many fields, including biological, technological, and social networks. Traditional stability methods often rely on explicit eigenvalue calculations or the construction of Lyapunov functions, which become computationally difficult for large interconnected systems. In this work, we introduce a stability classifier that provides a computationally efficient way to determine the stability of complex systems without performing full spectral analysis. Our approach combines the Dynamical Jacobian Ensemble framework, which uses the network structure and the underlying dynamics to construct a Jacobian, with the Gershgorin disc theorem, which provides bounds on the eigenvalues of the system's Jacobian matrix. Using this combination, we derive an analytical classifier that predicts whether a dynamical system is asymptotically stable, unstable, or sensitively stable based on a small set of dynamical and network parameters. The analysis reveals that the stability of large dynamical systems is strongly influenced by highly connected hub nodes. This observation allows the stability condition to be simplified by focusing on the contribution of the hub-associated Gershgorin disc. Extensive numerical simulations on large ensembles of Jacobian matrices confirm that the proposed classifier accurately separates stable, unstable, and sensitively stable regimes across a wide range of dynamical parameters and network structures. The resulting framework provides a simple and scalable tool for predicting the stability of large dynamical networks while capturing the interplay between network topology and interaction dynamics.
\end{abstract}

\section{Introduction}
Real-world systems are often large, complex, and continually exposed to perturbations that can propagate through their underlying interaction networks. In ecological communities \cite{CHEN_2024}, neural circuits \cite{SPORNS_2010,BASSETT_2017}, social systems \cite{DODDS_2005,CASTELLANO_2009}, and technological infrastructures \cite{DOBSON_2007,BULDYREV_2010},  perturbations may trigger regime shifts \cite{CHEN_2024,LENTON_2025}, cascading failures \cite{MOTTER_LAI_2002,CRUCITTI_2004,DOBSON_2007,BULDYREV_2010,DUAN_2019,SONG_2025}, or systemic collapse \cite{SCHEFFER_2009,BULDYREV_2010,CHEN_2024,LENTON_2025}. Understanding how complex systems respond to perturbations, and under what conditions they remain stable, is therefore a fundamental problem in dynamical systems theory, network science, and control theory \cite{SAUNDERS_1975,NEUBERT_1997,PASCUAL_2005,ZHI_2006,SCHOON_2012,KEFI_2019,BUMA_2020}.

Stability analysis provides a fundamental framework for characterizing the response of a system to perturbations. Different notions of stability can be distinguished according to the type of perturbation considered and the system’s response to it. For example, an equilibrium is Lyapunov stable if sufficiently small perturbations remain small, and asymptotically stable if they also decay to zero as time tends to infinity ~\cite{KHALIL_2002}. Classical approaches, including linear stability analysis~\cite{POINCARRE_1892}, Lyapunov theory~\cite{LYAPUNOV_ORIGINAL}, and the Routh-Hurwitz criterion~\cite{ROUTH_1877,HURWITZ_1895}, established a rigorous foundation for stability analysis, while structural stability theory \cite{PERRON_1929,GROBMAN_1959,HARTMAN_1960}, state-space methods \cite{KALMAN_1963}, stable manifold theory \cite{SMALE_1967}, and normal-form analysis \cite{ARNOLD_BOOK_1983} provide complementary frameworks for  stability analysis of nonlinear dynamical systems.

In the present work, we focus specifically on the local asymptotic stability of steady states. The stability of a steady state is assessed using the eigenvalues of the Jacobian matrix evaluated at that state. If all eigenvalues have strictly negative real parts, the steady state is locally asymptotically stable. 
Although conceptually straightforward, determining eigen spectrum becomes increasingly challenging for large networked systems, where the Jacobian can be high-dimensional and its spectral properties depend on both the underlying interaction topology and nonlinear dynamics.  This challenge, relevant to ecology \cite{MONTOYA_2006,SUWEIS_2015}, neuroscience \cite{SPORNS_2010,BASSETT_2017}, social systems \cite{CASTELLANO_2009,D_J_WATTS_BOOK_2003}, and technological networks \cite{BARABASI_1999,BARABASI_2016}, motivates scalable approaches for stability analysis in large networked systems. 

The relationship between complexity and stability in large dynamical systems has consequently attracted considerable attention for several decades. Gardner and Ashby first observed computationally that increasing system complexity can reduce stability \cite{GARDNER_ASHBY_1970}. May subsequently formalized this observation using random matrix theory (RMT), showing that sufficiently large and densely interacting random systems become increasingly unlikely to remain stable \cite{MAY_1972}. This seminal result established an influential theoretical connection between complexity and stability and provided a tractable framework for studying stability in high-dimensional systems. Subsequent studies, however, demonstrated that the stability of real-world systems cannot generally be understood from system size and interaction density alone. Network architecture, interaction patterns, correlations, and the underlying dynamical mechanisms can substantially alter stability \cite{SINHA_2005,ALLESINA_2012,TANG_PAWAR_ALLESINA_2014,KIRK_2015,OSULLIVAN_2019,TARNOWSKI_2020,BARBIER_2021}. These findings shifted the focus from purely random interactions toward understanding how the organization of network structure and interaction dynamics collectively governs stability.

In parallel, deterministic approaches such as decomposition-based Lyapunov methods \cite{SILJAK_1973,SILJAK_BOOK_1978} and Input-to-State Stability (ISS) \cite{SONTANG_1989} have extended stability analysis to large interconnected systems. More recent developments have addressed infinite-dimensional systems \cite{MIRONCHENKO_2020}, scalable network dynamics \cite{SILVA_2024}, and composite multiscale systems \cite{TANG_2026}. Despite these advances, the computational cost of spectral calculations and the system-specific construction of suitable Lyapunov functions can limit the practical applicability of conventional approaches to large complex networks. This motivates the search for alternative frameworks that can exploit the structural organization of networked dynamical systems to obtain scalable stability criteria.

A recent step in this direction was provided by Meena \textit{et al.} \cite{MEENA_2023} with the development of the Dynamical Jacobian Ensemble (DJE), which characterizes the Jacobian of a large networked dynamical system by incorporating information from both the network structure and the underlying interaction dynamics. The DJE captures the dependence of the Jacobian entries on these two components through a small set of dynamical and topological scaling exponents. In contrast to classical random-matrix approaches, this framework explicitly accounts for the dependence of the Jacobian entries on node connectivity and nonlinear interaction mechanisms, thereby establishing a direct connection between network structure, dynamics, and stability. Within the DJE framework, Meena \textit{et al.} \cite{MEENA_2023} derived a stability classifier based on the asymptotic behavior of the principal eigenvalue of the Jacobian, expressed in terms of the relevant dynamical and topological exponents. To obtain an analytical expression for the principal eigenvalue, their analysis exploited the strong degree heterogeneity of scale-free networks and employed a star-network approximation, in which the dominant hub capture the leading contribution in the large-network limit. While this approximation provides an analytically tractable description of the dominant role of the hub, it represents the original heterogeneous network through an idealized star topology. This motivates the development of an alternative approach for determining stability without relying on the star-network approximation.

Here, we combine the DJE \cite{MEENA_2023} with the Gershgorin Disc Theorem (GDT)  to develop an alternative framework for determining the local stability of large networked dynamical systems. The Gershgorin theorem provides bounds on the Jacobian eigenvalues directly from its matrix elements, thereby avoiding explicit computation of the complete eigenspectrum. For heterogeneous scale-free networks, we find that the Gershgorin disc associated with the most highly connected node frequently contains the principal eigenvalue, revealing the dominant contribution of the hub to stability. Based on this observation and the scaling properties of the DJE, we derive a Stability Classifier, $S$, in terms of the dynamical exponents $\Omega=(\eta,\mu,\nu,\rho)$ and topology-dependent scaling exponents. We validate the classifier using large ensembles of Jacobian matrices for different scale-free network structures and further examine its applicability to representative regulatory, population, mutualistic, inhibitory, epidemic, and biochemical dynamical systems \cite{MEENA_2023,BARZEL_2013}. Our framework thus provides a scalable approach for predicting the local stability of steady states in large networked dynamical systems without relying on the star-network approximation or explicit computation of the complete Jacobian spectrum.

\begin{figure*}[!h]
    \centering
    \includegraphics[width=0.95\linewidth]{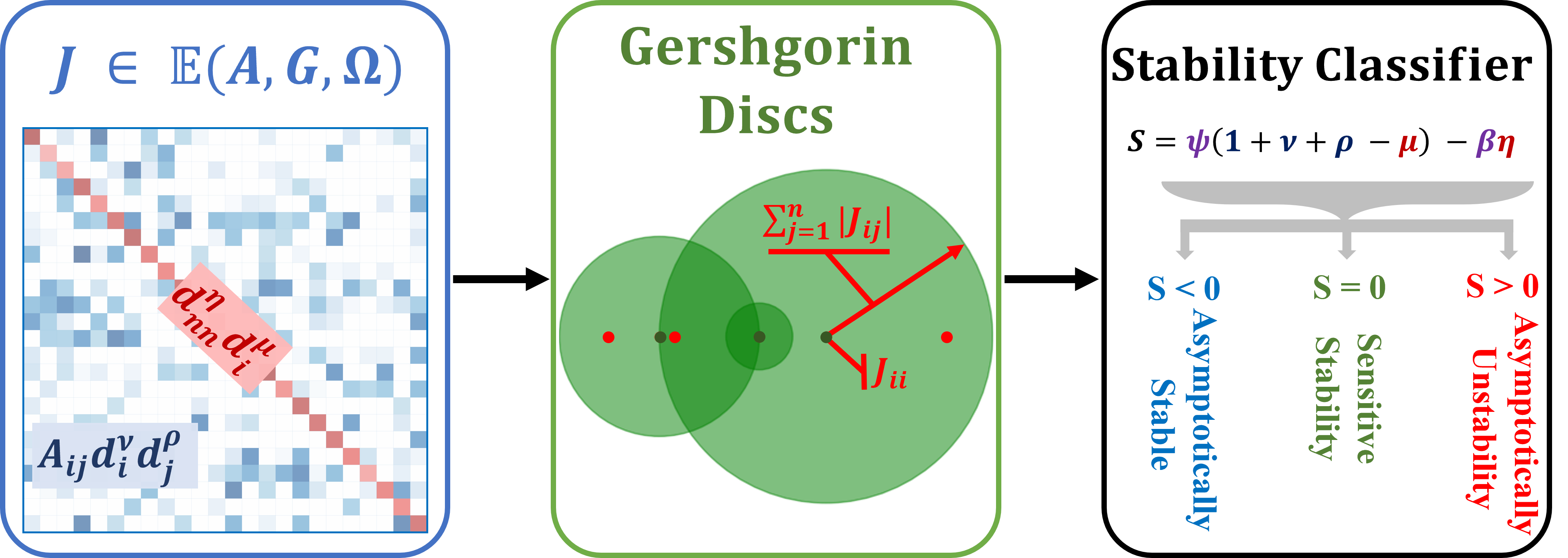}
\caption{
\textbf{Schematic illustration of the proposed framework combining the DJE and the Gershgorin Disc Theorem to construct the Stability Classifier.} \textbf{(a)} Many real-world complex systems can be modeled using the Barzel - Barabási equation, which provides a general dynamical framework for describing pairwise interactions in networked dynamical systems. Combining dynamics and network structure yields the Jacobian matrix $J \in E(A,G,\Omega)$ whose framework was laid in Meena et.al. \cite{MEENA_2023} . Unlike purely random ensembles, the entries of $J$ exhibit non-random scaling patterns. The non-zero elements correspond to existing network links, but their magnitudes depend on node degrees $d_i$ and $d_j$ as well as on the dynamical exponents ($\Omega =(\eta, \mu, \nu, \rho)$) \cite{MEENA_2023}. (\textbf{b}) The Gershgorin Disc Theorem is then applied to the Jacobian. Each disc is constructed using one row of $J$: the green dots represent the disc centers, which are the diagonal entries $J_{ii}$ of the Jacobian matrix, and the radii are given by the sum of the absolute values of the off-diagonal entries in that row.  These discs provide bounds on the eigenvalues which are shown as red dots without requiring full spectral computation. \textbf{(c)} Using these bounds, we derive the Stability Classifier $S$. The sign of $S$ determines the asymptotic behavior of the system: $S<0$ corresponds to asymptotic stability (blue), $S>0$ to asymptotic instability (red), and $S=0$ to sensitive stability (green). The classifier depends on the network heterogeneity exponents $\beta$ and $\psi$  together with the dynamical exponents $\Omega$ , providing a computationally efficient framework for stability prediction in large complex dynamical systems.
}
\label{fig:stability-diagram}
\end{figure*}

\section{Steady state Analysis}
\label{sec:FPA}

We consider the generalized framework for networked dynamical systems introduced by Barzel and Barabási~\cite{BARZEL_2013}, describing pairwise interactions among $N$ interacting components:

\begin{equation}
\frac{dx_i}{dt} = M_0(x_i(t)) + g \sum_{j=1}^{N} A_{ij}G_{ij} M_1(x_i(t)) M_2(x_j(t)).
\label{Master}
\end{equation}

Here $i = 1, \dots, N$, where $N$ is the number of nodes . The state of the system is represented by the vector $\mathbf{x}(t)=(x_1(t),\ldots,x_N(t))^T$, where $x_i(t)$ denotes the dynamical state of node $i$. The adjacency matrix $A$ encodes the network topology, with $A_{ij}=1$ if nodes $i$ and $j$ are connected and $A_{ij}=0$ otherwise, while $G_{ij}$ specifies the weight associated with the interaction between nodes $i$ and $j$. $g$ represents the global coupling strength. The function $M_0(x_i)$ describes the intrinsic dynamics of node $i$, such as growth or decay, whereas the product $M_1(x_i)M_2(x_j)$ captures the interaction mechanism between nodes $i$ and $j$. The functions $M_q(x)$, with $q=0,1,2$, share a common functional form across the network. Eq.~\eqref{Master} provides a general framework that encompasses a broad class of networked dynamical systems through appropriate choices of the functions $M_q(x)$. Examples include epidemic spreading models, where $x_i$ represents the infection probability~\cite{DODDS_2005,HUFNAGEL_2004,PASTOR_2001}; biochemical reaction dynamics, where $x_i$ represents the concentration of a reactant~\cite{VOIT_2000,MASLOV_2007,MASLOV_SERGEI_2007,YAN_MASLOV_2008}; birth--death processes, where $x_i$ represents the local population~\cite{GARDINER_2004,GANESH_2004,NOVOZHILOV_2006}; and gene-regulatory dynamics, where $x_i$ represents the expression level~\cite{ALON_BOOK_2006,KARLEBACH_2008}.

To determine the local asymptotic stability of Eq.~\eqref{Master}, we consider the dynamics in the vicinity of a steady state $\mathbf{x}_{\alpha}$. A steady state is a time-independent solution satisfying
 
\begin{equation}
\left.\frac{d\mathbf{x}}{dt}\right|_{\mathbf{x}=\mathbf{x}_{\alpha}}=0.
\end{equation}
Introducing a sufficiently small perturbation $\delta\mathbf{x}(t)$ around the steady state,
\begin{equation}
\mathbf{x}(t)=\mathbf{x}_{\alpha}+\delta\mathbf{x}(t),
\end{equation}
the dynamics can be linearized as
\begin{equation}
\frac{d\,\delta\mathbf{x}}{dt}
=
J\,\delta\mathbf{x}
+
\mathcal{O}\!\left(\|\delta\mathbf{x}\|^2\right),
\label{linearized}
\end{equation}
where $J$ is the $N\times N$ Jacobian matrix evaluated at $\mathbf{x}_{\alpha}$. Each element $J_{ij}$ quantifies the linear response of node $i$ to a small perturbation in node $j$. The local asymptotic stability of the steady state is determined by the eigenvalue spectrum of $J$. In particular, defining the principal eigenvalue $\lambda_{\max}$ as the eigenvalue with the largest real part, the steady state is locally asymptotically stable when \cite{STROGATZ_2018}
\begin{equation}
\operatorname{Re}(\lambda_{\max})<0.
\label{stability_condition}
\end{equation}

Conversely, $\operatorname{Re}(\lambda_{\max})>0$ implies local instability. Thus, the principal eigenvalue provides a direct spectral criterion for determining the local stability of the steady state. 
For dynamical networks, elements of Jacobian matrices interplay between network structures and dynamics \cite{MEENA_2023}, and mathematically expressed as
\begin{equation}
    J_{ij} =
    \begin{cases}
        -d_{nn}^{\,\eta}\, d_i^{\,\mu}, & \text{if } i = j, \\[6pt]
        d_i^{\,\nu}\, A_{ij}G_{ij}\, d_j^{\,\rho}, & \text{if } i \neq j,
    \end{cases}
    \label{Jacobian}
\end{equation}
where $d_i = \sum_{j=1}^{N} A_{ij}G_{ij}$ is the weighted degree of node $i$, and $d_{nn} = \frac{1}{N}\sum_{i=1}^{N}\frac{1}{d_i}\sum_{j=1}^{N}A_{ij}G_{ij}d_j$ is the average weighted nearest-neighbor degree. The dynamical exponents $\Omega=(\eta,\mu,\nu,\rho)$ depends on $M_q$ and determine how the Jacobian entries scale with network structural quantities, such as node degree and average nearest-neighbor degree \cite{MEENA_2023}.

\section{Graphical analysis of dynamical systems using Gershgorin discs}
\label{sec:GDT}
To analyze the local stability of the steady states, first we have to examine the eigenvalue spectrum of the Jacobian matrix. For large networked systems, direct computation of the complete eigenspectrum of an $N\times N$ Jacobian can become computationally demanding. We therefore employ the Gershgorin Disc Theorem (GDT), which provides bounds on the locations of the eigenvalues directly from the matrix elements without requiring explicit diagonalisation.

For a square matrix $M=[m_{ij}]$, the GDT states that every eigenvalue of $M$ lies within at least one Gershgorin disc (GD). The disc associated with row $i$ is defined as
\begin{equation}
D_i=
\left\{
\lambda\in\mathbb{C}:
|\lambda-m_{ii}|\leq R_i
\right\},
\label{gershgorin_disc}
\end{equation}
where $m_{ii}$ represents the center of the disc and
\begin{equation}
R_i = \sum_{j=1,\, j\neq i}^{N} |m_{ij}|.
\label{ri}
\end{equation}

is its radius. Consequently, the complete eigenspectrum of $M$, $\lambda(M)$, (refer figure \ref{fig:stability-diagram}(b)) satisfies

\begin{equation}
\lambda(M)\subseteq\bigcup_{i=1}^{N}D_i,
\label{gershgorin_union}
\end{equation}

When applied to the Jacobian $J$, the theorem therefore provides a direct geometrical bound on the locations of its eigenvalues in the complex plane. In particular, if all GDs lie entirely in the left half of the complex plane, all eigenvalues of $J$ necessarily have negative real parts, providing a sufficient condition for local asymptotic stability.

 \begin{table}[ht]
\centering

\begin{tabular}{c l c c}
\hline
$\gamma$ & Regime &
$\operatorname{Re}(\lambda_{\max})$ &
$\displaystyle\max_i\left(J_{ii}+\sum_{j=1}^{N}|J_{ij}|\right)$ \\
\hline
2.3 & Sensitively stable & $5.55\times10^{-17}$ & $0.17$ \\
2.3 & Stable            & $-0.6$& $-0.24$ \\
2.3 & Unstable          & $1.33\times10^{3}$ & $2.14\times10^{3}$ \\
\hline
2.5 & Sensitively stable & $-5.74\times10^{-16}$ & $0.22$ \\
2.5 & Stable            & $-0.20$ & $-0.07$ \\
2.5 & Unstable          & $11.19$ & $350.71$ \\
\hline
2.7 & Sensitively stable & $-4.15\times10^{-16}$ & $0.16$ \\
2.7 & Stable            & $-0.26$ & $-0.03$ \\
2.7 & Unstable          & $8.21$ & $132.47$ \\
\hline
\end{tabular}
\caption{Comparison between the maximum real part of the Jacobian eigenvalues, $\operatorname{Re}(\lambda_{\max})$, and the maximum bound put by the rightmost GD, $\max_i\left(J_{ii}+\sum_{j=1}^{N}|J_{ij}|\right)$, for different network degree exponents $\gamma$ and dynamical regimes. The stable cases exhibit negative values for both quantities, whereas the unstable cases exhibit positive values. For the sensitively stable cases, maximum eigenvalue real parts approximately equal to zero, and thus rightmost GD's maximum bound remains close to zero and positive.}
\label{tab:eigenvalue_gershgorin_comparison}
\end{table}

\begin{figure*}[!t]
    \centering
    \includegraphics[width=1.0\linewidth]{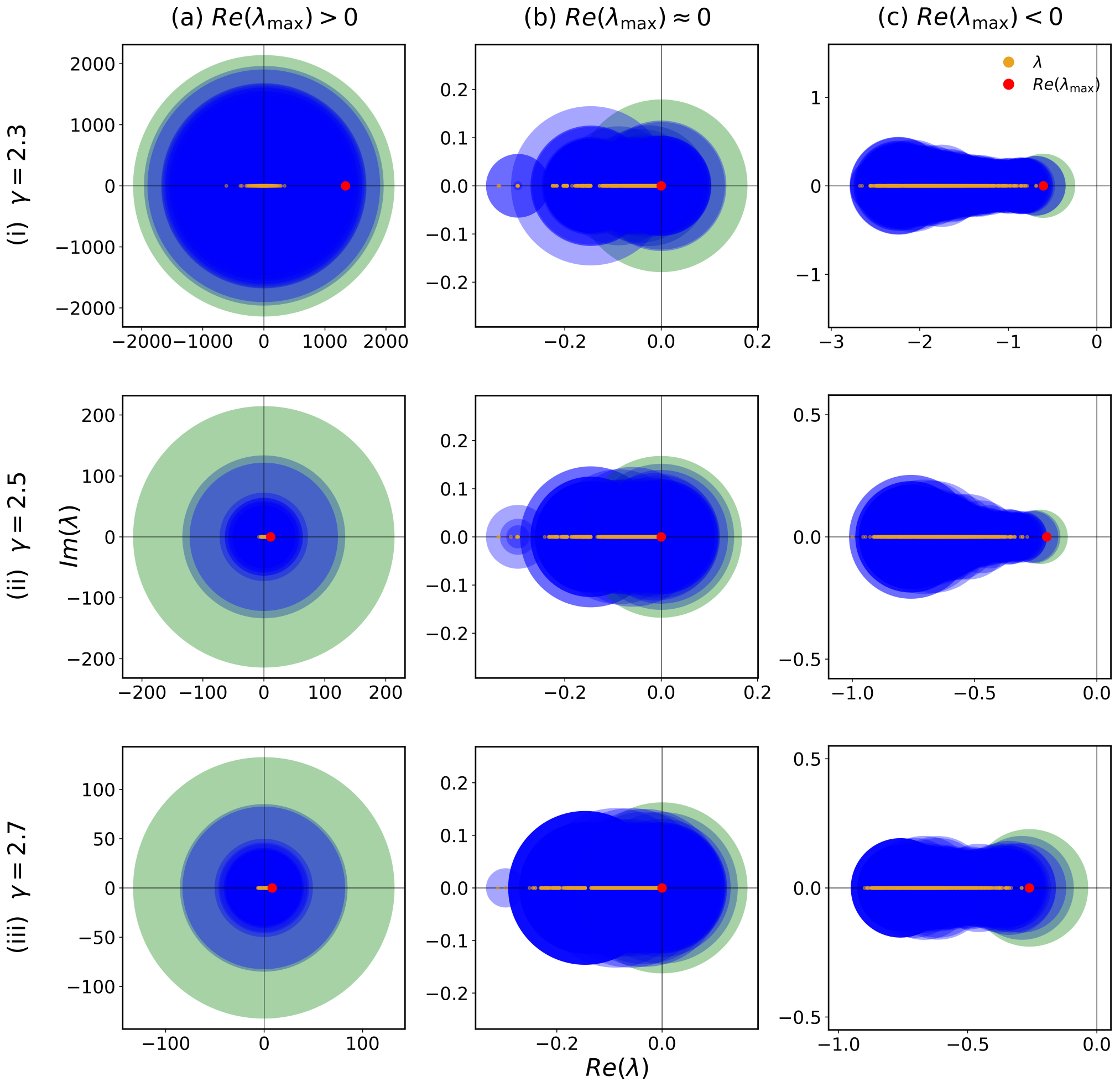}
    \caption{
    \textbf{Gershgorin discs for Dynamical Jacobian Ensembles associated with three stability regimes.} For these representative GDs, we consider three scale-free dynamical networks of with degree exponents $\gamma=2.3$, $2.5$, and $2.7$, corresponding to the first, second, and third rows, respectively. Here network size is $N=2000$ and average degree is $6$. The columns represent different dynamical stability regimes classified according to the principal eigenvalue $\lambda$ of the Jacobian matrix: the first column corresponds to the unstable regime ($\lambda>0$), the second column to the critical regime ($\lambda=0$), and the third column to the stable regime ($\lambda<0$). In the unstable regime, the GDs are centered close to the origin and extend across both sides of the imaginary axis. In the critical regime, some discs remain close to the origin, while others extend predominantly into the left half-plane. In the stable regime, the discs predominantly occupy the left half-plane, with greater coverage of the negative real half-plane than of the positive real half-plane. The GD associated with the hub node is highlighted in green, the remaining discs are shown in blue, and the principal eigenvalue is indicated by a red marker and the remaining spectral values are indicated by yellow color.}
    \label{fig:gershgorin}
\end{figure*}

\begin{figure*}[!t]
    \centering
    \includegraphics[width=1\linewidth]{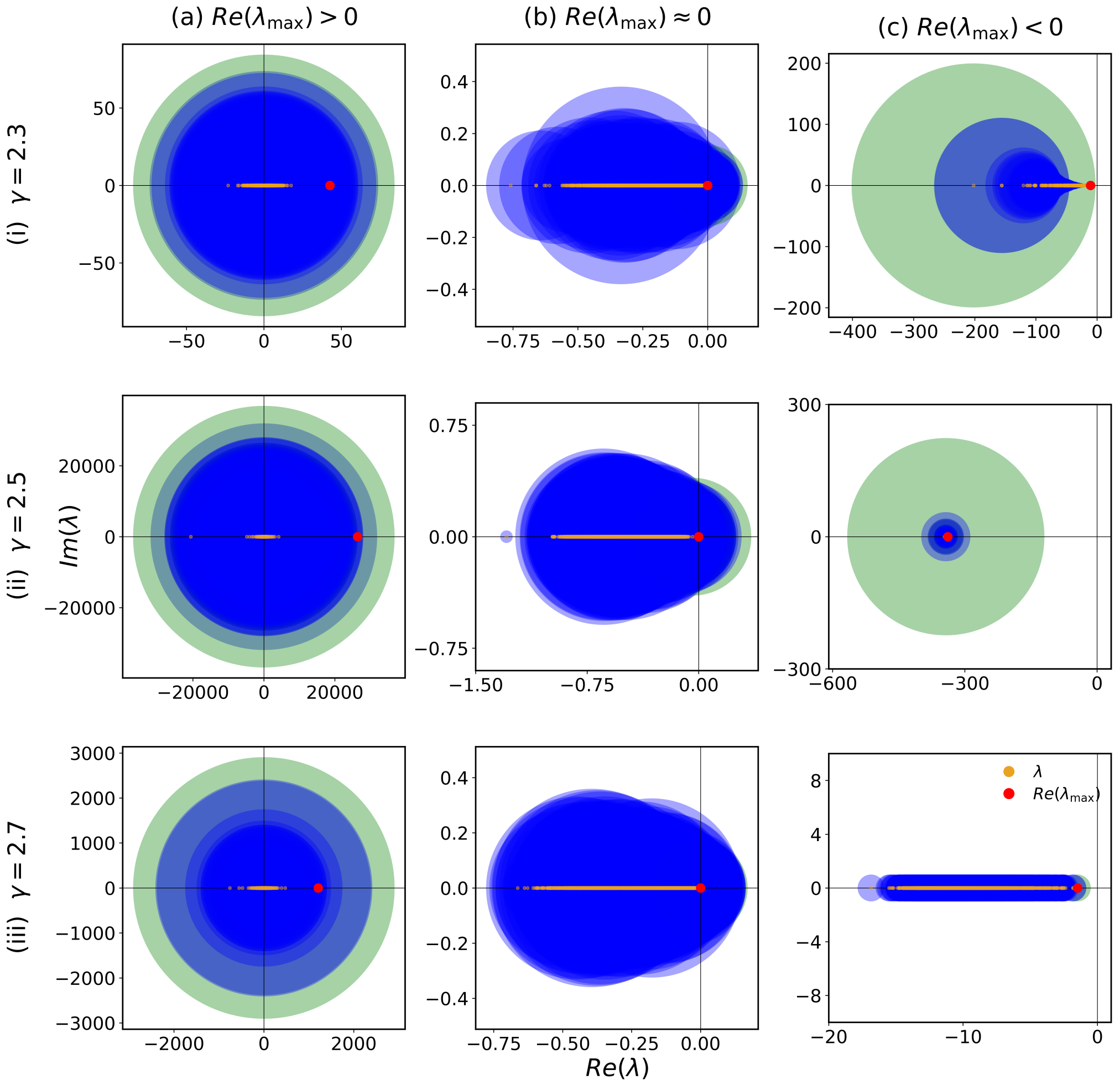}
\caption{Gershgorin discs for Dynamical Jacobian Ensembles associated with three stability regimes in weighted scale-free networks. The interaction weights are drawn from a normal distribution with mean 1 and standard deviation 0.1, i.e., \(G_{ij}\sim\mathcal{N}(1,0.1)\). The remaining network configurations, stability regimes, and graphical representations follow those described in Fig. \ref{fig:gershgorin}.}
    \label{fig:weighted_GD}
\end{figure*}

Our analysis focuses on scale-free networks of size $N=2000$ with  degree distribution $P(k)\propto k^{-\gamma}$, with $\gamma=2.3$, $2.5$, and $2.7$, as most empirical networks satisfy $2<\gamma<3$~ \cite{NEWMAN_BOOK_2018}. Primarily, we considered the global coupling constant, $g = 1$  and consider the network as an unweighted network, i.e., $G_{ij} = 1 \ \ \forall \ \ i,j$. For each $\gamma$, we vary the dynamical exponents $\Omega=(\eta,\mu,\nu,\rho)$, generating an ensemble of Jacobians which can be constructed using ~\eqref{Jacobian}. Further we extend our analysis to weighted networks as well by considering the weight distribution, $G_{ij} \approx  \mathcal{N}(1,0.1)$ (refer Fig.~\ref{fig:weighted_GD})

Some representative GDs corresponding to asymptotic stable, unstable and sensitive stability are  shown in Figs.~\ref{fig:gershgorin} and ~\ref{fig:weighted_GD}  that  reveals clear distinctions between stability regimes. 
Corresponding to the unstable regime ($\mathrm{Re}(\lambda_{\max}) > 0$), Figs.~\ref{fig:gershgorin} (a.i)–(a.iii) and ~\ref{fig:weighted_GD} (a.i)–(a.iii) show that the discs appear nearly symmetric about the imaginary axis, indicating that their centers are distributed around zero and that approximately equal area is covered on both sides. When $\mathrm{Re}(\lambda_{\max}) \approx 0$, the system enters a sensitively stable or near-critical regime (Figs. ~\ref{fig:gershgorin}(b.i)–(b.iii) and ~\ref{fig:weighted_GD}(b.i)–(b.iii)), where some discs remain centered near the origin while others extend into the left half-plane, producing a mixed spectral configuration. When $\mathrm{Re}(\lambda_{\max}) < 0$, corresponding to the stable regime (Fig.~\ref{fig:gershgorin}(c.i)–(c.iii) and ~\ref{fig:weighted_GD}(c.i)–(c.iii)), the union of GDs lies predominantly in the left half-plane, with significantly greater spectral coverage in the negative real half-plane than in the positive real half-plane. As shown in Table ~\ref{tab:eigenvalue_gershgorin_comparison}, the maximum rightmost GD's bound consistently captures the stable and unstable regimes across all three network degree exponents. For the sensitively stable cases, the bound lies nearby zero, consistent with the corresponding $Re(\lambda_{max})\approx 0$ . These geometric patterns provide a direct graphical separation between unstable, sensitive, and stable regimes.

We further observe in Figs.~\ref{fig:gershgorin} and ~\ref{fig:weighted_GD} that, across the stable, unstable, and near-critical regimes, the dominant eigenvalue is located within the GD associated with the hub node shown in green color which is also evident in the case of weighted network This correspondence is observed consistently for the three scale-free networks considered, with $\gamma=2.3$, $2.5$, and $2.7$. Thus, the position of the hub-associated disc closely follows the location of the dominant eigenvalue across the different stability regimes, indicating that this disc provides a bound on the principal eigenvalue and will help to determine the stability of system. Therefore, we investigate when hub-associated bounds capture the dominant spectral behavior and assess their accuracy and limitations across network structures and dynamical regimes.

 \begin{figure}
     \centering
     \includegraphics[width=1\linewidth]{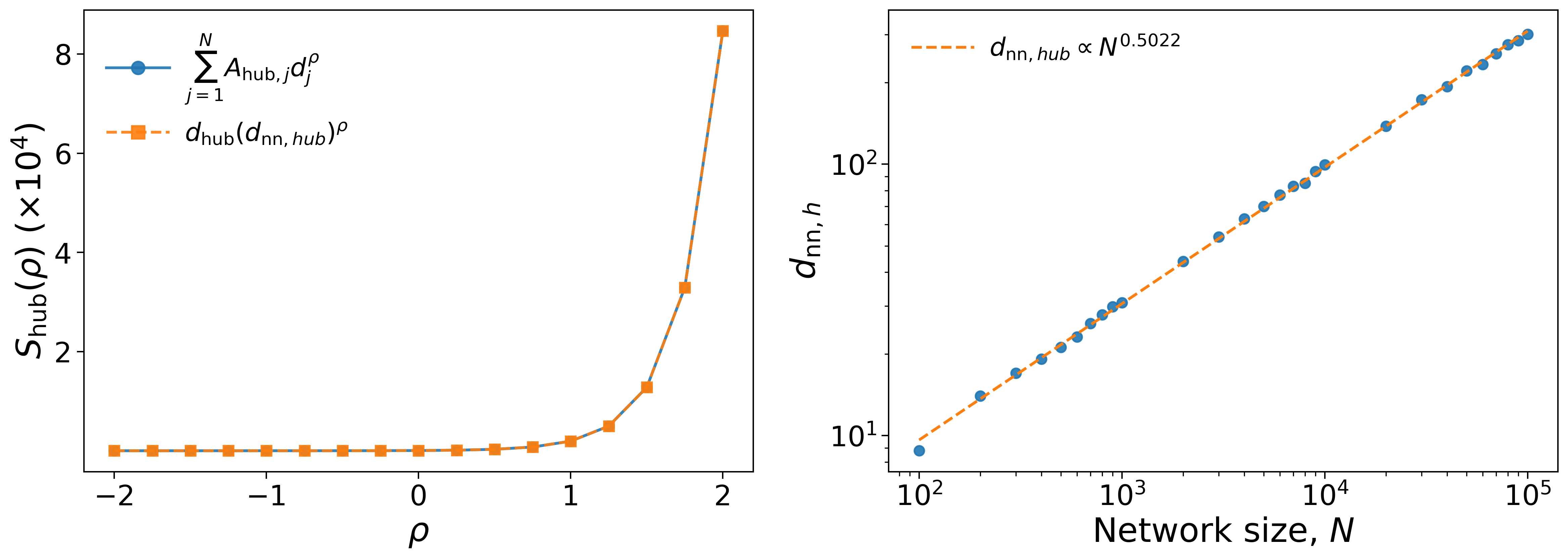}
     \caption{\textbf{Approximation of the hub-based moment \(S_{\mathrm{hub}}(\rho)\) and scaling of the mean degree of the hub's nearest neighbors with network size.} (a) Comparison between the actual moment, \(\sum_{j=1}^{N} A_{\mathrm{hub},j}d_j^\rho\) (blue circles), and its approximation, \(d_{\mathrm{hub}}(d_{\mathrm{nn,hub}})^\rho\) (orange squares), as a function of the exponent \(\rho\). The close overlap between the two curves indicates that the approximation accurately captures the hub-based moment over the investigated range of \(\rho\). (b) The log–log plot of the mean degree of the hub's nearest neighbors, \(d_{\mathrm{nn,hub}}\), against network size \(N\). The dashed orange line represents the power-law scaling \(d_{\mathrm{nn,hub}}\propto N^{0.5022}\) which yields the scaling exponent $\psi_{nn} = 0.5022$.}
     \label{d_nn,hub_scaling}
 \end{figure}
 
\section{Dynamical Stability Classifier}
\label{DSC}
Using the DJE framework and the GDT, introduced in sections ~\ref{sec:FPA} and \ref{sec:GDT}, respectively, 
we construct a stability classifier, denoted by $S$, to characterize the 
asymptotic stability of large-scale networked dynamical systems. The 
classifier combines the scaling properties of the network topology with 
the dynamical exponents governing the entries of the Jacobian matrix, 
thereby providing a computationally efficient alternative to the explicit 
spectral decomposition of large Jacobian matrices.  

Rather than explicitly computing 
$Re(\lambda_{max)}$, we employ the GDT to constrain 
its location in the complex plane.

For the $i$th row of $J$, the corresponding GD is centered at 
$J_{ii}$ and has radius
\begin{equation}
    R_i=\sum_{j=1,\,j\neq i}^{N}|J_{ij}|.
\end{equation}
A sufficient condition to lie all eigenvalues within the 
left half of the complex plane is
\begin{equation}
    J_{ii}+R_i<0.
    \label{GDT_stability_general}
\end{equation}
Since the diagonal entries of the DJE considered 
here are negative, Eq.~\eqref{GDT_stability_general} can equivalently be 
written as
\begin{equation}
    |J_{ii}|>R_i.
    \label{GDT_stability}
\end{equation}

For a dynamical J, the diagonal and off-diagonal 
contributions are expressed in terms of the network degrees and dynamical 
exponents\cite{MEENA_2023}. Substitution into Eq.~\eqref{GDT_stability} gives
\begin{equation}
    d_{nn}^{\eta}d_i^{\mu}
    >
    \sum_{j=1,\,j\neq i}^{N}
    A_{ij}G_{ij}d_i^{\nu}d_j^{\rho},
    \label{inequality}
\end{equation}
where $d_i$ is the degree of node $i$, $d_{nn}$ is the average 
nearest-neighbor degree, and $\eta$, $\mu$, $\nu$, and $\rho$ are the 
dynamical exponents characterizing the corresponding dynamical system. 

The Dynamical Jacobian Ensembles considered 
here exhibit an important structural feature. As shown in 
Fig.~\ref{fig:gershgorin}, the principal eigenvalue is consistently found within 
the GD associated with the hub node. This behavior is observed 
for different scale-free degree exponents and across the all dynamical stability regimes considered in our analysis. 
Thus, 

the hub-associated GD provides a useful bound for tracking the 
location of the principal eigenvalue

Motivated by this systematic observation, we construct the stability 
classifier from the Gershgorin condition associated with the hub node. 
Setting $i=\mathrm{hub}$ in Eq.~\eqref{inequality} gives
\begin{equation}
    d_{nn}^{\eta}d_{\mathrm{hub}}^{\mu-\nu}
    >
    \sum_{j=1}^{N}
    A_{\mathrm{hub},j}|G_{\mathrm{hub},j}|d_j^{\rho}.
    \label{hub_condition}
\end{equation}

For an uncorrelated network, we can use the following approximation:
\begin{equation}
    \sum_{j=1}^{N}
    A_{\mathrm{hub},j}|G_{\mathrm{hub},j}|d_j^{\rho}
    \approx d_{\text{hub}}d_{nn,hub}^{\rho}
    \label{approx}
\end{equation}
where,  
\begin{equation}
    d_{nn,h}
    =
    \frac{1}{d_{\mathrm{hub}}}
    \sum_{j=1}^{N}
    A_{\mathrm{hub},j}d_j.
    \label{hub_neighbor_degree}
\end{equation}
Thus, eq ~\eqref{hub_condition} can be written as, 
\begin{equation}
     d_{nn}^{\eta}d_{\mathrm{hub}}^{\mu-\nu} > d_{\mathrm{hub}}d_{nn,\mathrm{hub}}^\rho
     \label{uncorrelated_network_aprrox}
\end{equation}
Fig.~\ref{d_nn,hub_scaling}(a) verifies the approximation given by Eq.~\eqref{hub_neighbor_degree}, while Fig.~\ref{d_nn,hub_scaling}(b) demonstrates the scaling relation \(d_{nn,\mathrm{hub}}\sim N^{\psi_{nn}}\) is obtained by examining how the hub's average nearest-neighbor degree varies with network size with \(\psi_{nn}\) denoting the corresponding scaling exponent. 
Furthermore, for the scale-free networks considered here, the hub degree and the average nearest-neighbor degree scale as
\begin{equation}
    d_{\mathrm{hub}}\sim N^{\psi},
    \label{psi_scaling}
\end{equation}
and
\begin{equation}
    d_{nn}\sim N^{\beta}
    \label{beta_scaling}
\end{equation}

For highly heterogeneous networks, we find that $\psi_{nn} \approx \psi$. As $\gamma$ increases, the network tends to be more homogeneous and the degree distribution of the network also changes which disobeys  $\psi_{nn} \approx \psi$ . In such scenario, $\psi_{nn}$ can be calculated independently using Eq.~\eqref{hub_neighbor_degree} and checking how does it scales with the network size.
Substituting Eqs.~\eqref{psi_scaling}-\eqref{beta_scaling} and re-arranging eq~\eqref{uncorrelated_network_aprrox} yields, 
\begin{equation}
       N^{\beta\eta  + \psi (\mu - \nu - \rho - 1 )}>1.
\end{equation}
Thus, for a system to be stable, $\psi(1 + \nu + \rho-\mu) - \beta\eta$, needs to be less than zero. This motivates the definition of the stability classifier, $S$,
 \begin{equation}
     S = \psi(1 + \nu + \rho-\mu) - \beta\eta
     \label{Stability_Classifier}
 \end{equation}

The parameters \(\psi\) and \(\beta\) are key network scaling exponents that describe the growth behavior of the hub degree and the average nearest-neighbor degree with increasing system size \(N\). 
The stability classifier defined in Eq.~\eqref{Stability_Classifier} has been tested on the scale free network. The three regimes of stability are thus clearly distinct based on the following conditions. 
\begin{equation}
S  
\begin{cases}
< 0 & \text{(Stable)} \\
= 0 & \text{(Sensitively stable)} \\
> 0 & \text{(Unstable)} 
\end{cases}
\end{equation}
Here, the dynamical exponents $\Omega = (\eta, \mu, \nu, \rho)$ are determined by the underlying interaction functions $M_q(x)$, $q=0,1,2$, while the topological scaling exponents $\Omega' = (\beta, \psi)$ depend only on the network structure and are estimated from degree moments (refer sections~\ref{sec:beta} and~\ref{sec:psi}). For Eq.~\eqref{Stability_Classifier}, we employed a hub-node approximation by focusing on the GD associated with the most highly connected node. As shown in Section~\ref{sec:GDT}, the dominant eigenvalue $\lambda_{\max}$ lies within the hub-associated disc in most cases. Motivated by this observation, we evaluated the Gershgorin stability condition using only the Jacobian row corresponding to the hub node, rather than all rows. This approximation simplifies the analysis without significantly reducing its predictive accuracy, and leads naturally to the classifier $S$. Eq.~\eqref{Stability_Classifier} therefore provides a compact analytical framework that integrates dynamical parameters and network topology. Thus we can use $S$ to predict a priori whether a specific combination of topology and dynamics will exhibit stable functionality or not.

\begin{figure*}[!t]
    \centering
    \includegraphics[width=1\linewidth]{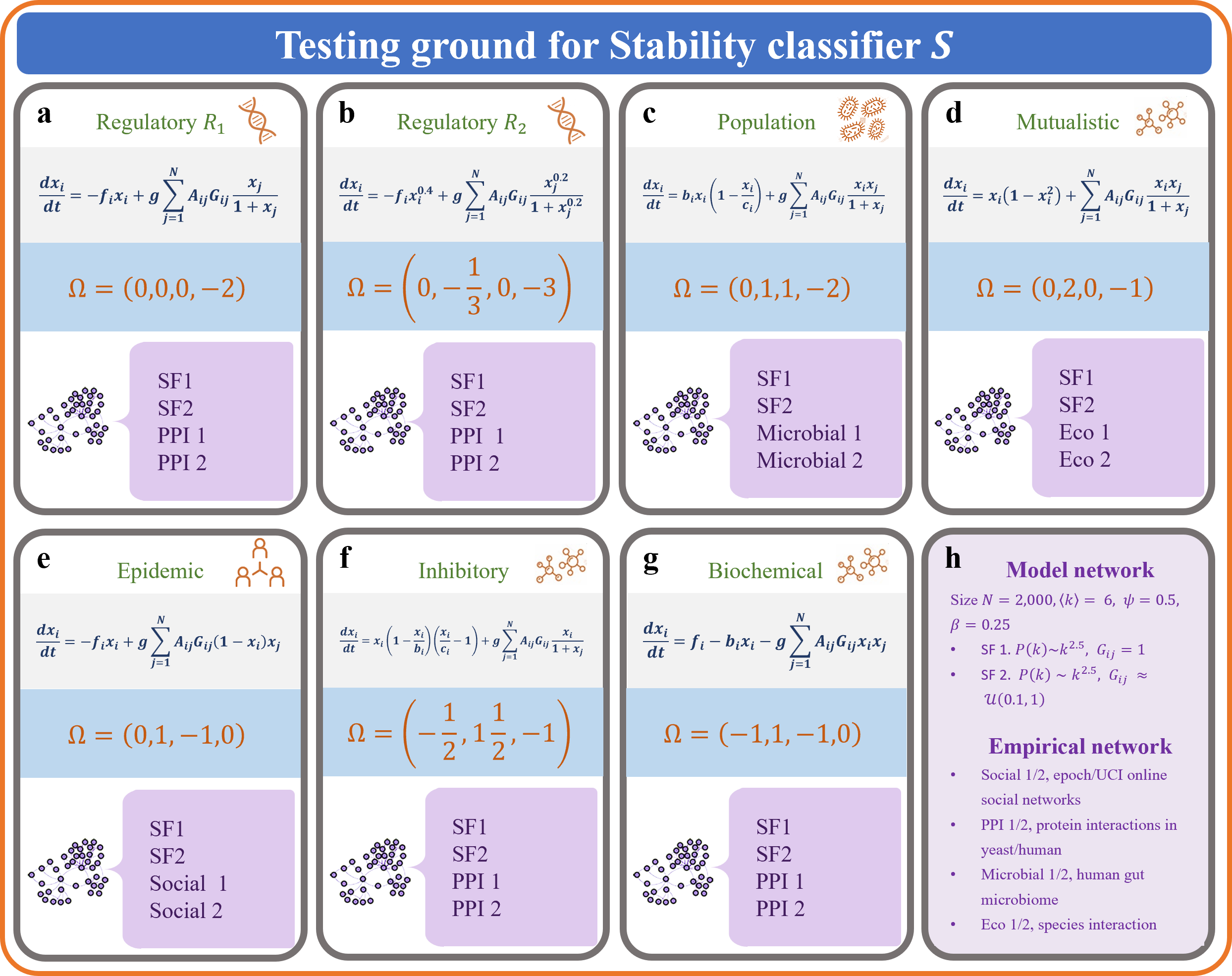}
    \caption{\textbf{Testing ground for stability classifier, $S$.} We tested our stability classifier for different dynamical systems using the Dynamical exponents derived using the DJE framework. \textbf{a,b, Regulatory:} Subcellular dynamics following the Michaelis--Menten model. This system exhibits an active fixed points in which all $x_i>0$ (namely, pandemic). In this active state, the system has a set of exponents $\Omega$ derived using the Dynamical Jacobian Ensemble. Here we present $\Omega$ and $S$ for the active fixed point. The remaining panels follow a similar format. Here we have two cases of the same dynamical model, by changing the exponents of interacting state variables, $x_i $ and $x_j$, to show that how the intrinsic parameter affects the dynamics and hence the stability. \textbf{c, Population}: Population dynamics, for example in microbial communities. \textbf{d, Mutualistic}: Mutualistic interaction dynamics, where species benefit from positive interactions with others, for example plant–pollinator systems or cooperative microbial communities. \textbf{e, Inhibitory:} Suppression dynamics, for example between hosts and pathogens. \textbf{f, Epidemic:} Susceptible-infected-susceptible (SIS) dynamics \textbf{g, Biochemical:} Protein--protein interactions modeled via mass-action kinetics. This system exhibits a single fixed point. \textbf{h, Networks:} We considered a modeled scale free network of size 2000 with an average degree of 6 and tested these above dynamical systems}
    \label{fig:Real}
\end{figure*}

\section{Stability Classifier on Real world complex systems}
 
To examine the predictive stability of our proposed stability classifier $S$  on real world complex systems, we constructed a broad testing ground, including relevant dynamic models from different domains. Epidemic, the susceptible–infected susceptible model \cite{PASTOR_2015, DODDS_2003, BROCKMANN_2009} for disease spreading; Regulatory, the Michaelis–Menten model \cite{KARLEBACH_2008} for gene regulation; Inhibitory, growth suppression in pathogen–host interaction \cite{WODARZ_2002}; Biochemical, protein–protein interaction \cite{MURRAY_1989,BARUCH_2011,BARUCH_2012} in sub- cellular networks; Population, a model of species \cite{HOLLING_1959} interaction in population dynamics ; and Mutualistic, to track the abundance of species in a mutualistic environment \cite{HOLLAND_2002}. For each dynamical system we have a unique set of dynamical exponents. These dynamical exponents are derived by following the algorithm that has been extensively discussed in Meena \textit{et al.}  \cite{MEENA_2023} corresponding to the Eq.~\eqref{Master}. Figure \ref{fig:Real} mentions unique $\Omega$ for every dynamical system.

We used an uncorrelated, unweighted, undirected scale-free network with size $N=2000$, average degree $\langle k\rangle = 6$, and  $\gamma = 2.5$. We consider an uncorrelated, undirected scale-free network of size $N = 2000$, with an average degree $\langle k\rangle = 6$, and $\gamma = 2.5$. We examine two configurations, an unweighted network, $G_{ij} = 1$, and a weighted network, where the interaction weights are drawn from a uniform distribution, i.e., i.e.,  $G_{ij} \sim \mathcal{U}(0.1, 1.0)$. The initial conditions for all nodes are independently drawn from a uniform distribution in the interval $[0.1, 1.0]$, i.e.,  $x_i(0) \sim \mathcal{U}(0.1, 1.0) \ \ \forall \ \ i$ $\epsilon$ $N$. Using these initial conditions, the system reached a steady state after following initial transients. After perturbing the system around a fixed point, the Jacobian of all the dynamical systems mentioned in the Fig \ref{fig:Real} is constructed by linearizing the dynamical equations around their steady states.  The principal eigenvalue, $Re(\lambda_{max}$) is calculated for Jacobian corresponding to each dynamical system and it's stability was assessed using Eq.~\eqref{stability_condition} and compared with $S$. The proposed Stability Classifier \(S\) demonstrates high accuracy in predicting the stability of dynamical systems across both unweighted and weighted network configurations. Its consistent performance across the two interaction-weight settings supports the reliability and applicability of the framework for stability analysis in complex networks.

\begin{figure*}[!t]
    \centering
    \includegraphics[width=1\linewidth]{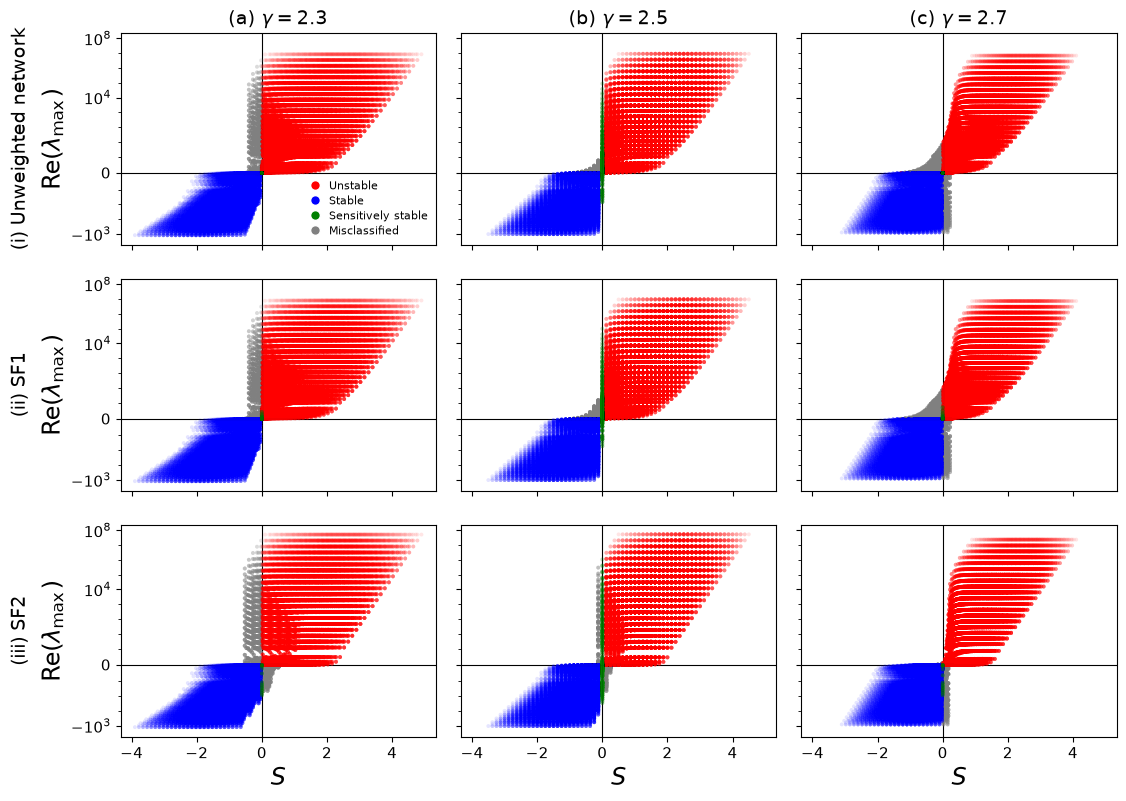}
 \caption{\textbf{Three classes of dynamic stability.}
Relationship between the Stability Classifier $S$ and the principal eigenvalue $\mathrm{Re}(\lambda)$ for $83,521$ Jacobian samples for each network configuration. The Jacobian matrices are generated by varying the dynamical exponents $\eta, \mu, \nu, \rho$ across a broad region of parameter space. The underlying interaction topology is a scale-free network with an average degree of $6$ and a network size of $2,000$ nodes. Three network configurations are considered: \textbf{i)} an unweighted network ( $G_{ij}=1$); \textbf{ii)} SF1, in which the interaction weights are drawn from a normal distribution, $G_{ij}\sim\mathcal{N}(1,0.1)$; and \textbf{iii)} SF2, in which the interaction weights are drawn from a power-law distribution, $P(G)\sim G^{-3}$. For each configuration, the three columns correspond to \textbf{a)} $\gamma=2.3$, \textbf{b)} $\gamma=2.5$, and \textbf{c)} $\gamma=2.7$. $S$ classifies the DJE into three distinct regimes: unstable ($S>0$; red), stable ($S<0$; blue), sensitively stable ($S=0$; green). We also identified misclassified cases shown in gray.}
\label{fig:Re_vs_S}
\end{figure*}

\section{Stability Classifier on Dynamical Jacobian Ensemble}
The Stability classifier $S$ in ~\eqref{Stability_Classifier}  helps group all the \( J \in \mathbb{E}(A, G, \Omega) \) into distinct classes by identifying the relevant topological, $\Omega' = (\beta , \psi)$ and dynamical, $\Omega  = (\eta, \mu, \nu, \rho)$  control parameters that help analytically predict the stability of any system within the form of the ~\eqref{Master} . 
To test the predictive power of  $S$, we tested it on uncorrelated scale-free network with \(N=2000\) with $\psi = 0.5$, average degree \(\langle k\rangle=6\), and degree exponents \(\gamma=2.3,\ 2.5,\) and \(2.7\) considering three network configurations: an unweighted network with \(G_{ij}=1\), a weighted network, \(G_{ij}\sim\mathcal{N}(1,0.1)\), and a weighted network with interaction weights following a power-law distribution, \(P(G)\sim G^{-3}\). In particular, using DJE, 83,561 Jacobians are generated for each type of network configuration keeping the network structure constant and calculated the corresponding principal eigenvalue, $Re(\lambda_{max})$ for each Jacobian. Simultaneously, we also calculate the stability classifier, $S$ using the ~\eqref{Stability_Classifier} for specific $\Omega  = (\eta, \mu, \nu, \rho)$ and $\Omega' = (\beta , \psi)$ and compared it with the corresponding $Re(\lambda_{max})$. As predicted graphically, we find that \( S \) sharply splits the sample into three classes. The asymptotically unstable class (figure \ref{fig:Re_vs_S}, red, top right) has \( S > 0 \) and consequently also \( \operatorname{Re}(\lambda) > 0 \), indicating guaranteed instability. The asymptotically stable class (figure \ref{fig:Re_vs_S}, blue, bottom left) is observed for \( S < 0 \), and has, in all cases, \( \operatorname{Re}(\lambda) < 0 \), corresponding to stable dynamics. Finally, for \( S = 0 \), we observe sensitive stability, with \( \operatorname{Re}(\lambda) \) showing no clear asymptotic behavior, either positive or negative (figure \ref{fig:Re_vs_S}, green). Table ~\ref{tab:accuracy} display the accuracy achieved for different values of $\gamma$ and weight distributions. A percentage (\(6.8\%\) for unweighted networks, \(7.17\%\) for networks with normal weight distribution, and \(4.77\%\) for power-law weighted networks, yielding an average error of \(6.25\%\)) remains misclassified due to the hub-node approximation that we considered while deriving the expression for the stability classifier, as while analyzing graphically, there were some cases when $Re(\lambda_{max})$ didn't lie in the hub node associated GD.  The number of misclassified cases decreases as we consider a highly heterogeneous network  (refer Fig.~\ref{fig:Re_vs_S}).  Although a small fraction of cases did not follow the hub-node approximation, these exceptions highlight the richness and diversity of behaviors that can arise in complex networks. 

\begin{figure*}[!t]
    \centering
    \includegraphics[width=1\linewidth]{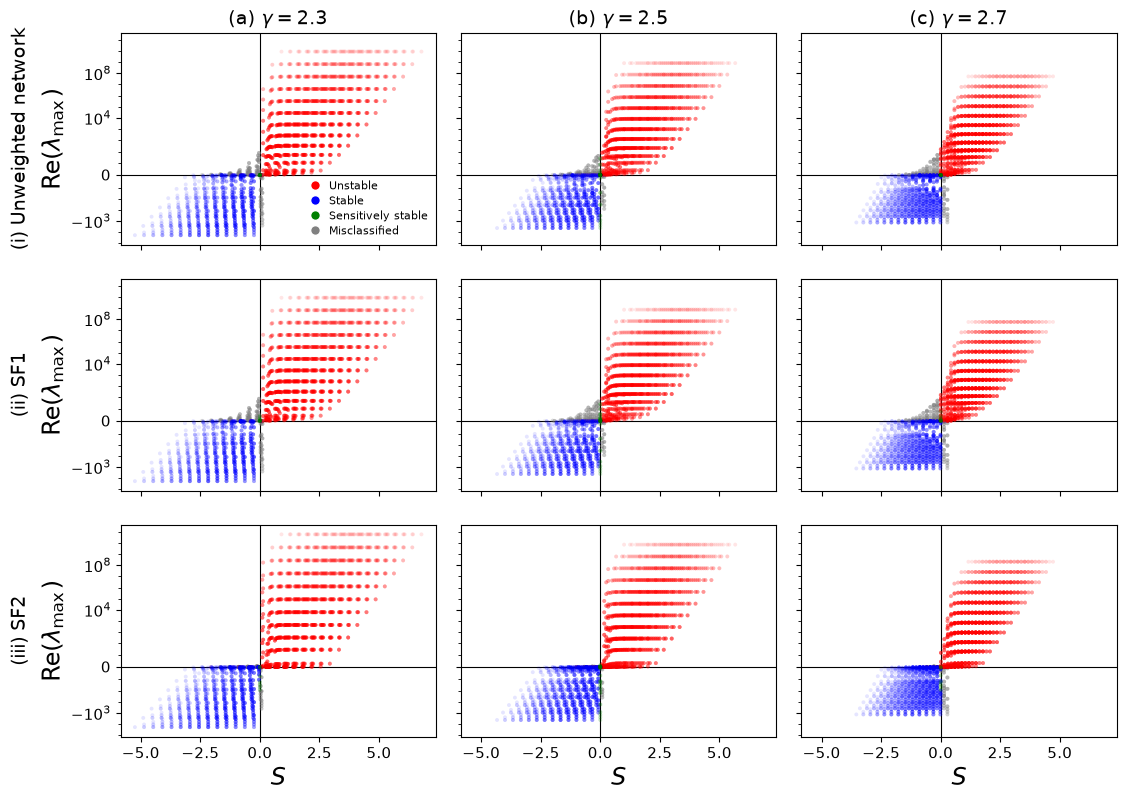}
 \caption{\textbf{Three classes of dynamic stability.}
Relationship between the Stability Classifier $S$ and the principal eigenvalue $\mathrm{Re}(\lambda)$ for 6,561 Jacobian samples for each network configuration. The Jacobian matrices are generated by varying the dynamical exponents $\eta, \mu, \nu, \rho$ across a broad region of parameter space. The underlying interaction topology is a scale-free network with an average degree of $6$ and a network size of $2,000$ nodes. For a correlated network, three network configurations are considered: \textbf{i)} an unweighted network ( $G_{ij}=1$); \textbf{ii)} SF1, in which the interaction weights are drawn from a normal distribution, $G_{ij}\sim\mathcal{N}(1,0.1)$; and \textbf{iii)} SF2, in which the interaction weights are drawn from a power-law distribution, $P(G)\sim G^{-3}$. For each configuration, the three columns correspond to \textbf{a)} $\gamma=2.3$, \textbf{b)} $\gamma=2.5$, and \textbf{c)} $\gamma=2.7$. $S$ classifies the DJE into three distinct regimes: unstable ($S>0$; red), stable ($S<0$; blue), sensitively stable ($S=0$; green). We also identified misclassified cases shown in gray.}
\label{corr_S}
\end{figure*}

To examine the applicability of this approximation beyond uncorrelated networks, we also test $S$ on correlated networks of size $N = 2000$ following $d_{\text{hub}} \sim N^{\frac{1}{\gamma -1}}$, average degree \(\langle k\rangle=6\), and degree exponents \(\gamma=2.3,\ 2.5,\) and \(2.7\) considering three network configurations: an unweighted network with \(G_{ij}=1\), a weighted network, \(G_{ij}\sim\mathcal{N}(0.1,1.0)\), and a weighted network with interaction weights following a power-law distribution, \(P(G)\sim G^{-3}\). We achieve an accuracy of $\sim92.9\%$ in the case of correlated networks for all type of network configuration mentioned above (refer Fig. ~\ref{corr_S}). This verifies that our Stability classifier is not restricted for uncorrelated networks only. 

\begin{table}[ht]
    \centering
    \begin{tabular}{cccc}
        \hline
        $\gamma$&
        $G_{ij}=1$&
        $G_{ij}\sim\mathcal{N}(0.1,1.0)$&
        $P(G)\sim G^{-3}$ \\
        \hline
        2.3&  94.81& 94.45& 96.54\\
        2.5&  93.16& 93.12& 96.69\\
        2.7&  91.62& 90.91& 92.46\\
        \hline
    \end{tabular}
    \caption{Classification accuracy of the Stability Classifier $S$ for scale-free networks with different degree exponents $\gamma$ and interaction-weight configurations. The entries report the percentage of correctly classified dynamical systems for unweighted networks ($G_{ij}=1$) and weighted networks with interaction weights drawn from a uniform distribution, $G_{ij}\sim\mathcal{N}(0.1,1.0)$, or a power-law distribution, $P(G)\sim G^{-3}$.}
    \label{tab:accuracy}
\end{table}

\begin{figure*}[!t]
    \centering
    \includegraphics[width=1.0\linewidth]{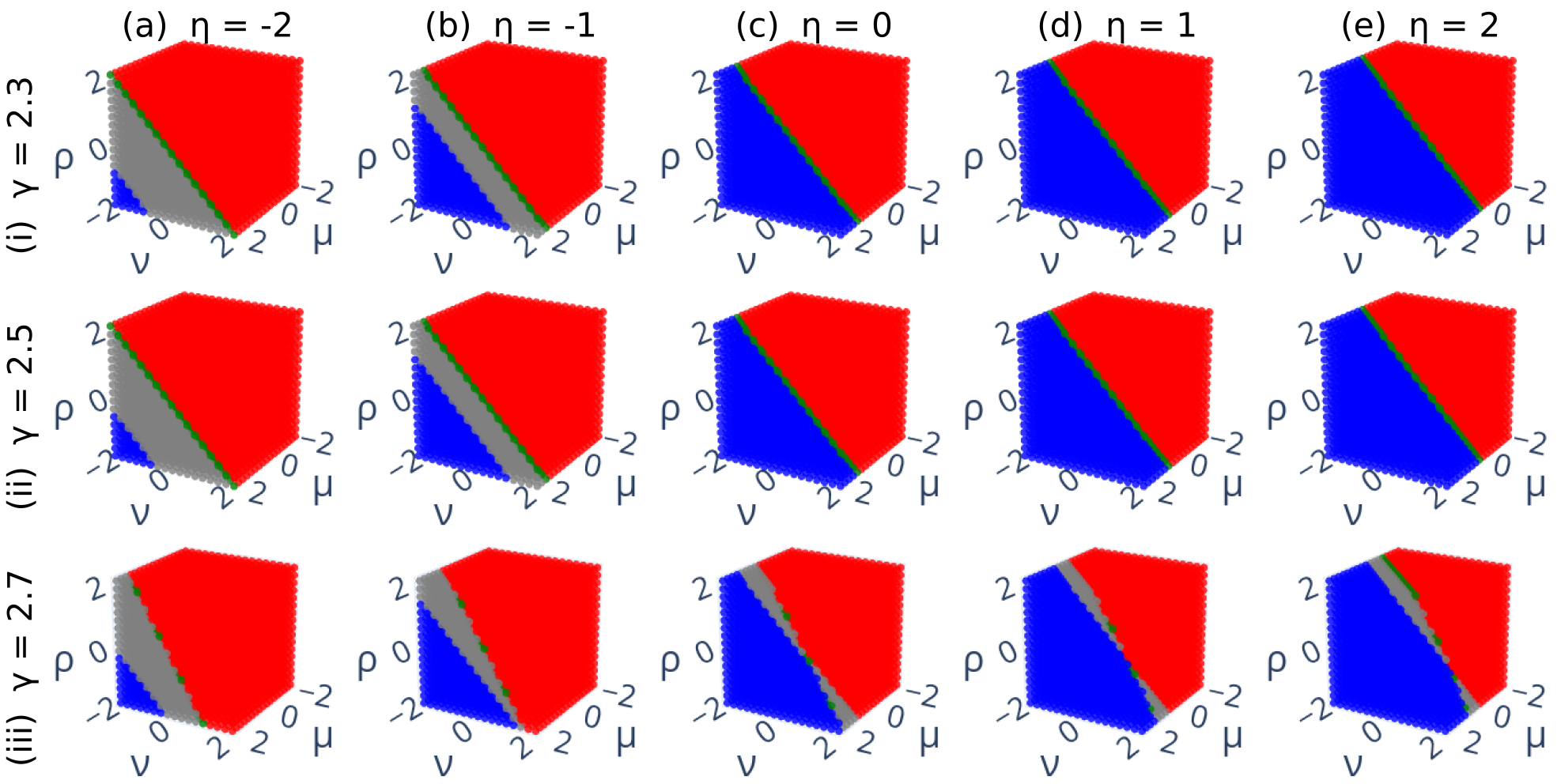}
\caption{
Parameter-space representation of the Stability Classifier $S$ for scale-free networks with degree exponents (i) $\gamma=2.3$, (ii) $\gamma=2.5$, and (iii) $\gamma=2.7$. Columns correspond to different values of the dynamical exponent $\eta$: (a) $\eta=-2$, (b) $\eta=-1$, (c) $\eta=0$, (d) $\eta=1$, and (e) $\eta=2$. For each combination of $(\gamma,\eta)$, the parameters $\mu$, $\nu$, and $\rho$ are varied to explore the stability landscape in the $(\mu,\nu,\rho)$ space. Blue regions correspond to asymptotically stable systems ($S<0$), red regions to asymptotically unstable systems ($S>0$), and the green surface ($S=0$) represents the sensitively stable or marginal stability boundary separating the two regimes. Gray regions indicate parameter combinations for which the Stability Classifier misclassifies the dynamical system.
}
    \label{3D_plots}
\end{figure*}

\subsection*{Stability Classifier in Parameter Space}

To understand how dynamical and topological exponents affect stability, we explore the parameter space of $S$. By varying $\mu$, $\nu$, and $\rho$ while keeping $\eta$ and the network topology fixed, i.e., $\Omega'=(\beta,\psi)$ fixed, stability maps is generated that show how dynamical and topological exponents influence system's behavior.We use the same values of \(S\) and \(Re(\lambda_{\max})\) obtained for Fig. \ref{fig:Re_vs_S} to generate the parameter space.As shown in Fig.~\ref{3D_plots}, the sign of $S$ divides the $(\mu,\nu,\rho)$ space into distinct regions: $S<0$ correspond to asymptotically stable dynamics (blue), $S>0$ correspond to unstable dynamics (red), and $S=0$ defines a marginal or sensitively stable boundary (green) separating the two regimes. Transitions between these regimes occur whenever $S$ changes sign. The gray regions denote parameter combinations for which the Stability Classifier misclassifies the dynamical system.

Fig.~\ref{3D_plots} reveals two clear trends. First, for a fixed network topology, increasing $\eta$ progressively reduces both the gray and red regions, while the blue stable region expands and eventually dominates the parameter space. This trend is evident in Fig.~\ref{3D_plots} (i)((a)-(e)), corresponding to $\gamma=2.3$, and similar behavior is observed in rows (ii) and (iii), corresponding to $\gamma=2.5$ and $\gamma=2.7$, respectively. The reduction of the gray region indicates an improvement in the predictive capability of the Stability Classifier, whereas the shrinking of the unstable region suggests that larger values of $\eta$ favor stability. This behavior can be understood from the role of $\eta$, which controls the diagonal entries of the Jacobian matrix. Since the eigenvalue spectrum is strongly influenced by the diagonal terms, increasing $\eta$ enhances the contribution of self-dynamics relative to the off-diagonal interaction terms. Consequently, the GD centers move further into the stable region, increasing the likelihood that the dominant eigenvalue has a negative real part and thereby stabilizing the system. This highlights the importance of self-dynamics in governing stability. A second trend emerges as the degree exponent $\gamma$ increases. This effect is reflected in Fig.~\ref{3D_plots} (c-e) ((i)-(iii)), where the gray misclassified region decreases as $\gamma$ decreases, indicating an improvement in the prediction accuracy of the Stability Classifier for more heterogeneous networks for $\eta \geq 0 $ cases where the vice versa trend can be seen in Fig.~\ref{3D_plots} (a-b) ((i)-(iii)) indicating that the diagonal entries affects the stability of the dynamical system . At the same time,  Fig.~\ref{3D_plots} (b) ((i)-(iii)), (c) ((i)-(iii)), (d) ((i)-(iii)), and (e) ((i)-(iii)) show a systematic increase in the blue stable region. Overall, the parameter-space representation demonstrates that the Stability Classifier provides a direct link between internal dynamical mechanisms and large-scale network stability.

\section{Network Parameters}
\label{sec:stability_classifier_parameters}

The stability classifier $S$ in Eq.~\eqref{Stability_Classifier}
reduces the dynamical stability of Eq.~\eqref{Master} to six
relevant exponents: the dynamical exponents
$\Omega=(\eta,\mu,\nu,\rho)$ and the topological exponents
$\Omega'=(\beta,\psi)$. The first four exponents,
$\Omega=(\eta,\mu,\nu,\rho)$, are determined by the system's intrinsic
dynamics $M_q(x)$ around each of its fixed points~\cite{MEENA_2023}.
Thus, $\Omega$ is embedded in the natural interaction dynamics of the
system, as determined by the governing equations and the manner in
which the nodes influence one another through the dynamical model.
The remaining two exponents, $\Omega'=(\beta,\psi)$, characterize the
network topology and are determined by the scaling properties of the
adjacency matrix $A$ and the weight matrix $G$.

The exponent $\beta$ characterizes the scaling of the average
nearest-neighbour degree $d_{nn}$ with network size $N$, whereas
$\psi$ characterizes the scaling of the degree of the hub node,
$d_{\mathrm{hub}}$. These exponents quantify how degree heterogeneity
and hub dominance evolve as the network grows. Their estimation
therefore provides a means of relating network structure to the
stability characteristics of large-scale dynamical systems.

\subsection{\texorpdfstring{Degree Heterogeneity exponent, $\beta$}Degree Heterogeneity exponent, $\beta$}
\label{sec:beta}
The degree heterogeneity exponent $\beta$ quantifies the dependence of the average nearest-neighbour degree $d_{nn}$ on the network size $N$. We express this dependence through the power-law scaling
\begin{equation}
    d_{nn}\sim N^\beta.
    \label{beta_scaling_supp}
\end{equation}
This scaling
relation provides a measure of how the connectivity of a typical
neighbour changes as the network grows.

For an uncorrelated network, the average nearest-neighbour degree is
related to the first and second moments of the degree distribution
$P(d)$ through~\cite{GAO_2016,YAN_2017}
\begin{equation}
    d_{nn}=\frac{\langle d^2\rangle}{\langle d\rangle},
    \label{eq:dnn_base}
\end{equation}
where $\langle d\rangle$ and $\langle d^2\rangle$ denote the first
and second moments of the degree distribution, respectively.

For an uncorrelated scale-free network with degree distribution
$P(d)\sim d^{-\gamma}$, where $2<\gamma<3$, the second moment depends
on the lower and upper degree cutoffs. Denoting the minimum degree
by $d_{\min}$ and the degree of the hub by $d_{\mathrm{hub}}$, we
obtain
\begin{equation}
    \langle d^2\rangle
    =\int_{d_{\min}}^{d_{\mathrm{hub}}}d^2P(d)\,dd
    \sim\int_{d_{\min}}^{d_{\mathrm{hub}}}d^{2-\gamma}\,dd.
    \label{eq:second_moment_integral}
\end{equation}
Evaluating the integral gives
\begin{equation}
    \begin{aligned}
        \langle d^2\rangle
        &\sim
        \left[
        \frac{d^{3-\gamma}}{3-\gamma}
        \right]_{d_{\min}}^{d_{\mathrm{hub}}}\\
        &=
        \frac{d_{\mathrm{hub}}^{3-\gamma}
        -d_{\min}^{3-\gamma}}{3-\gamma}.
    \end{aligned}
    \label{eq:second_moment_integrated}
\end{equation}
Here, $d_{\min}$ is assumed to remain fixed as the network size
$N$ increases. Since $2<\gamma<3$, the exponent $3-\gamma$ is
positive. Consequently, for sufficiently large $d_{\mathrm{hub}}$,
the upper-cutoff contribution dominates, yielding
\begin{equation}
    \langle d^2\rangle\sim d_{\mathrm{hub}}^{3-\gamma}.
    \label{eq:second_moment_hub}
\end{equation}

Assuming that the hub degree follows the scaling relation
$d_{\mathrm{hub}}\sim N^\psi$, Eq.~\eqref{eq:second_moment_hub}
becomes
\begin{equation}
    \langle d^2\rangle\sim N^{\psi(3-\gamma)}.
    \label{eq:second_moment_scaling}
\end{equation}
If the mean degree $\langle d\rangle$ remains finite and
approximately independent of $N$, substituting
Eq.~\eqref{eq:second_moment_scaling} into
Eq.~\eqref{eq:dnn_base} gives
\begin{equation}
    d_{nn}\sim N^{\psi(3-\gamma)}.
    \label{eq:dnn_scaling}
\end{equation}

Defining $\beta'=3-\gamma$, we obtain
\begin{equation}
    d_{nn}\sim N^{\psi\beta'}=N^\beta,
    \label{eq:dnn_beta}
\end{equation}
which yields the theoretical relation
\begin{equation}
    \boxed{\beta=\psi(3-\gamma)=\psi\beta'.}
    \label{eq:beta_theory}
\end{equation}
Thus, under the stated assumptions, $\beta$ quantifies the scaling
of the average nearest-neighbour degree with network size. For a
fixed value of $\psi$, smaller values of $\gamma$ correspond to
larger values of $\beta$, indicating a stronger dependence of
$d_{nn}$ on $N$.

The degree exponent $\gamma$ can be estimated from the slope of the
degree distribution on a log-log scale, while $\psi$ can be estimated
from the scaling of the hub degree with network size. Eqs.
\eqref{eq:dnn_scaling} and \eqref{eq:beta_theory} then provide the
corresponding theoretical prediction for $\beta$ under the stated
assumptions. Finite network sizes and differences in network
generation procedures may lead to deviations from this asymptotic
prediction. In particular, the relation
$\beta=\psi(3-\gamma)$ is conditional on the assumed power-law
degree distribution, the dominance of the upper degree cutoff, and
the approximate independence of the mean degree from $N$.

\subsection{\texorpdfstring{Hub scaling exponent, $\psi$}{Hub scaling exponent, psi}}
\label{sec:psi}

The exponent $\psi$ characterizes how the degree of the hub node,
$d_{\mathrm{hub}}$, scales with network size $N$. In scale-free
networks, hub nodes correspond to the most highly connected vertices
and play a dominant role in shaping both the structural and dynamical
properties of the system~\cite{BARABASI_1999}. Since the stability
classifier is derived using a hub-node approximation, understanding
the scaling behaviour of the hub degree is essential. We express
this dependence through the power-law relation
\begin{equation}
    d_{\mathrm{hub}}\sim N^\psi,
    \label{psi_scaling_s}
\end{equation}
where $\psi$ quantifies the growth rate of the hub degree with
network size on a logarithmic scale.

For a structural cutoff in an uncorrelated network, the maximum
degree scales as $d_{\mathrm{hub}}\sim N^{1/2}$, corresponding to
$\psi=0.5$. This cutoff constrains the largest degree to suppress
degree correlations that may otherwise arise when hubs become
excessively connected. For scale-free networks with $2<\gamma<3$,
the natural cutoff, which arises from the finite network size,
generally scales differently from the structural cutoff. The
appropriate scaling therefore depends on the network-generation
procedure and the constraints imposed on the degree distribution.

\begin{figure}[h]
    \centering
    \includegraphics[width=1\linewidth]{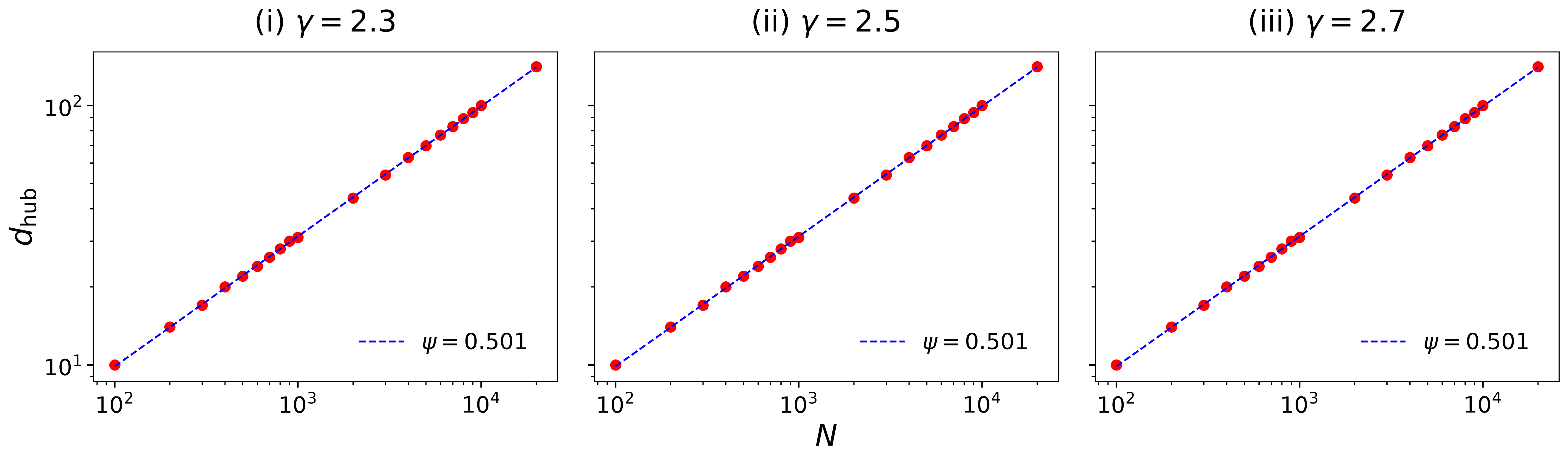}
    \caption{Scaling of the hub-node degree,
    $d_{\mathrm{hub}}$, with network size $N$ for uncorrelated
    scale-free networks with degree exponents $\gamma=2.3$,
    $2.5$, and $2.7$. All three networks exhibit the same
    empirically measured scaling exponent, $\psi=0.5$.
    This value is close to the structural-cutoff prediction
    $\psi=1/2$ and characterizes the growth of the hub degree
    with network size under the network configurations considered.}
    \label{fig:hub_node_scaling}
\end{figure}

Fig.~\ref{fig:hub_node_scaling} illustrates the measured dependence
of $d_{\mathrm{hub}}$ on $N$ for the uncorrelated scale-free networks
considered in this study.

Our numerical results for the three degree exponents,
$\gamma=2.3$, $2.5$, and $2.7$, exhibit a common empirical scaling
exponent $\psi=0.501$. This value is close to the structural-cutoff
prediction $\psi=0.5$, indicating that the hub degree grows
approximately as the square root of the network size in the
considered network configurations. The observed scaling is
consistent with the imposed structural-cutoff scaling, although
the measured exponent is empirical and may depend on the finite
network sizes and network-generation procedure.

The exponent $\psi$ therefore characterizes the observed dependence
of the hub degree on network size. Within the stability-classifier
framework, it serves as a topological parameter describing how the
hub contribution to the Jacobian changes as the network grows.
Together with $\beta$, it captures complementary aspects of network
structure: $\beta$ describes the scaling of the average
nearest-neighbour degree, while $\psi$ describes the scaling of
the most highly connected node.

For Erd\H{o}s--R\'enyi (ER) networks, the degree distribution is
binomial and approaches a Poisson distribution in the large-$N$
limit when the connection probability is held fixed. In this
setting, the average nearest-neighbour degree does not exhibit
the systematic power-law scaling with network size that
characterizes the scale-free networks considered above.
Consequently, the degree heterogeneity exponent is taken to be
$\beta=0$ within the power-law scaling description adopted here.

Similarly, ER networks do not exhibit the scale-free hub dominance
considered in the hub-node approximation. Their maximum degree
does not follow the same power-law scaling as the hubs of the
scale-free networks studied here. Accordingly, the hub scaling
exponent is taken to be $\psi=0$ within the stability-classifier
framework for ER networks. This assignment describes the absence
of power-law hub-degree scaling; it does not imply that the
maximum degree of a finite ER network is strictly independent
of $N$.

Substituting $\beta=0$ and $\psi=0$ into
Eq.~\eqref{Stability_Classifier}, all three cases reduce to
$S=0$. This corresponds to the sensitively stable regime, in
which stability neither systematically increases nor decreases
with network size and instead depends on the intrinsic rate
parameters of the large dynamical system. This result is
consistent with the findings reported in~\cite{MEENA_2023},
where ER networks are shown to exhibit sensitive stability
for large network sizes.

By incorporating the dynamical exponents
$\Omega=(\eta,\mu,\nu,\rho)$, the stability classifier accounts
for variations in the intrinsic driving mechanisms of the
dynamical system. The topological exponents
$\Omega'=(\beta,\psi)$ describe how network heterogeneity and
hub dominance change with network size. Their combined
contribution enables the stability classifier to capture the
interplay between the internal dynamics and the underlying
network structure. Therefore, $S$ incorporates both dynamical and topological
information to characterize the stability class of the system
around a specific fixed point: stable, unstable, or sensitively
stable. The estimation of $\beta$ and $\psi$ provides the
network-scaling information required to apply the classifier
to large-scale networks and to compare stability across
different network structures.

\section{Discussion}

In this work, we introduced a Stability Classifier that provides a unified and computationally efficient framework for analyzing the stability of complex networked dynamical systems. By combining the DJE with the GDT, the method avoids full spectral computation while still giving accurate predictions.  The dynamical exponents $\Omega = (\eta, \mu, \nu, \rho)$, derived from the interaction functions, together with the topology-dependent scaling exponents $\Omega' = (\beta, \psi)$, allow stability to be expressed in terms of both structure and dynamics. This formulation enables a systematic understanding of how internal interaction mechanisms and network topology jointly determine stability.

An important distinction between the present work and the earlier stability classifier proposed by Meena \textit{et al.}~ \cite{MEENA_2023} lies in the treatment of network topology. In the earlier framework, the stability classifier was derived analytically using a star-network approximation, which becomes valid in the thermodynamic limit where the degree of the dominant hub grows sufficiently large. While this approximation captures the leading contribution of highly connected nodes, it replaces the original scale-free structure by an idealized star topology. In contrast, the present work directly considers networks generated from a scale-free degree distribution and derives the stability classifier without invoking the star-network approximation. Consequently, the resulting framework preserves the heterogeneous structure of scale-free networks throughout the analysis and provides a more general description of the interplay between network topology and nonlinear dynamics.

A central result of this study is the dominant role of hub nodes in scale-free networks. The GD associated with the hub node often contains the dominant eigenvalue, allowing stability to be inferred from the behavior of this single row of the Jacobian. It also highlights the impact of degree heterogeneity: in fat-tailed networks, highly connected nodes strongly influence asymptotic stability. This explains why, in such systems, stability properties can become largely independent of microscopic parameters and instead be governed by topological scaling laws~ \cite{MEENA_2023}. A related phenomenon appears in epidemic dynamics on scale-free networks, where disease transmission and the basic reproduction number $R_0$ are strongly influenced by network topology and the presence of highly connected nodes ~\cite{ZHANG_2024}.  

Despite its ability to accurately classify stability across large complex systems, the framework has limitations. The hub-based approximation assumes that stability is primarily determined by the most connected node. While this assumption holds for most of the cases studied here, the Stability Classifier also exhibits misclassifications for a subset of dynamical systems, indicating that the stability of these systems is not determined solely by the hub. In such cases, other structural or dynamical features can contribute to the location of the dominant eigenvalue and therefore influence stability. These misclassifications provide an important indication that, although the hub plays a major role in determining stability in most scale-free networks, additional mechanisms can also contribute meaningfully to the stability of the system. Incorporating contributions from multiple influential nodes or relaxing the hub-based approximation could therefore improve the accuracy of the classifier in such cases.

The proposed framework has broad implications for real-world systems, including biological, technological, ecological, and engineered networks. By identifying how structural and dynamical factors influence stability, the classifier provides insight into which nodes promote stabilization or destabilization. This understanding can guide strategies such as protecting critical nodes or modifying interaction strengths to maintain stable operation. Moreover, the computational efficiency of the method makes it suitable for large-scale systems where full spectral analysis is impractical. Thus, the Stability Classifier establishes a direct connection between network topology and nonlinear dynamics for systems described by Eq.~\eqref{Master}. By reducing high-dimensional spectral analysis to a tractable analytical condition, it provides a practical and conceptually clear framework for studying stability in large complex dynamical systems, particularly those with scale-free or fat-tailed structures.

\bibliographystyle{unsrt}
\bibliography{ref}

@article{BARUCH_2011,
  title = {Binomial Moment Equations for Stochastic Reaction Systems},
  author = {Barzel, Baruch and Biham, Ofer},
  journal = {Phys. Rev. Lett.},
  volume = {106},
  issue = {15},
  pages = {150602},
  numpages = {4},
  year = {2011},
  month = {Apr},
  publisher = {American Physical Society},
  doi = {10.1103/PhysRevLett.106.150602},
  url = {https://link.aps.org/doi/10.1103/PhysRevLett.106.150602}
}

@book{SMALE_1967,
  title={Differentiable Dynamical Systems},
  author={Smale, Stephen},
  journal={Bulletin of the American Mathematical Society},
  volume={73},
  number={6},
  pages={747--817},
  year={1967}
}

@book{ARNOLD_BOOK_1983,
  title={Geometrical Methods in the Theory of Ordinary Differential Equations},
  author={Arnold, V. I.},
  publisher={Springer},
  year={1983}
}

@article{CHEN_2024,
title = {Stability of ecological systems: A theoretical review},
journal = {Physics Reports},
volume = {1088},
pages = {1-41},
year = {2024},
note = {Stability of Ecological Systems: A Theoretical Review},
issn = {0370-1573},
doi = {https://doi.org/10.1016/j.physrep.2024.08.001},
author = {Can Chen and Xu-Wen Wang and Yang-Yu Liu},
}

@article{CRUCITTI_2004,
  title = {Model for cascading failures in complex networks},
  author = {Crucitti, Paolo and Latora, Vito and Marchiori, Massimo},
  journal = {Phys. Rev. E},
  volume = {69},
  issue = {4},
  pages = {045104(R)},
  numpages = {4},
  year = {2004},
  month = {Apr},
  publisher = {American Physical Society},
  doi = {10.1103/PhysRevE.69.045104},
  url = {https://link.aps.org/doi/10.1103/PhysRevE.69.045104}
}

@article{MOTTER_LAI_2002,
  title = {Cascade-based attacks on complex networks},
  author = {Motter, Adilson E. and Lai, Ying-Cheng},
  journal = {Phys. Rev. E},
  volume = {66},
  issue = {6},
  pages = {065102(R)},
  numpages = {4},
  year = {2002},
  month = {Dec},
  publisher = {American Physical Society},
  doi = {10.1103/PhysRevE.66.065102},
  url = {https://link.aps.org/doi/10.1103/PhysRevE.66.065102}
}

@article{
DUAN_2019,
author = {Dongli Duan  and Changchun Lv  and Shubin Si  and Zhen Wang  and Daqing Li  and Jianxi Gao  and Shlomo Havlin  and H. Eugene Stanley  and Stefano Boccaletti },
title = {Universal behavior of cascading failures in interdependent networks},
journal = {Proceedings of the National Academy of Sciences},
volume = {116},
number = {45},
pages = {22452-22457},
year = {2019},
doi = {10.1073/pnas.1904421116},
URL = {https://www.pnas.org/doi/abs/10.1073/pnas.1904421116},
}

@article{DOBSON_2007,
    author = {Dobson, Ian and Carreras, Benjamin A. and Lynch, Vickie E. and Newman, David E.},
    title = {Complex systems analysis of series of blackouts: Cascading failure, critical points, and self-organization},
    journal = {Chaos: An Interdisciplinary Journal of Nonlinear Science},
    volume = {17},
    number = {2},
    pages = {026103},
    year = {2007},
    month = {06},
    issn = {1054-1500},
    doi = {10.1063/1.2737822},
    url = {https://doi.org/10.1063/1.2737822},

}

@article{BULDYREV_2010,
  author = {Buldyrev, Sergey V. and Parshani, Roni and Paul, Gerald and Stanley, H. Eugene and Havlin, Shlomo},
  title = {Catastrophic cascade of failures in interdependent networks},
  journal = {Nature},
  year = {2010},
  volume = {464},
  number = {7291},
  pages = {1025--1028},
  doi = {10.1038/nature08932},
  url = {https://doi.org/10.1038/nature08932},
  issn = {1476-4687},
  month = 04
}

@article{SONG_2025,
title = {Cascading failure model and resilience-based sequential recovery strategy for complex networks},
journal = {Reliability Engineering \&
System Safety},
volume = {253},
pages = {110488},
year = {2025},
issn = {0951-8320},
doi = {https://doi.org/10.1016/j.ress.2024.110488},
url = {https://www.sciencedirect.com/science/article/pii/S095183202400560X},
author = {Xudan Song and Pengcheng Zhao and Rongrong Yin and Yunxiao Zu and Yong Zhang},
}

@misc{LENTON_2025,
address = {Exeter, UK},
publisher = {University of Exeter},
title = {The Global Tipping Points Report 2025},
url = {https://pure.iiasa.ac.at/id/eprint/20918/},
author = {Lenton, T. M. and Milkoreit, M. and Wilcock, S. and Abrams, J. F. and Armstrong McKay, D. I. and Buxton, J. E. and Donges, J. F. and Loriani, S. and Wunderling, N. and Alkemade, F. and Barrett, M. and Constantino, S. and Powell, T. and Smith, S. R. and Boulton, C. A. and Pinho, P. and Dijkstra, H. A. and Pearce-Kelly, P. and Roman-Cuesta, R. M. and Dennis, D.}
}

@article{BASSETT_2017,
  author  = {Bassett, Danielle S. and Sporns, Olaf},
  title   = {Network neuroscience},
  journal = {Nature Neuroscience},
  year    = {2017},
  volume  = {20},
  number  = {3},
  pages   = {353--364},
  month   = mar,
  doi     = {10.1038/nn.4502},
  url     = {https://doi.org/10.1038/nn.4502},
  issn    = {1546-1726}
}

@book{SPORNS_2010,
author = {Sporns, Olaf},
year = {2010},
publisher = {MIT Press},
month = {01},
pages = {},
title = {Networks of the Brain}
}

@article{CASTELLANO_2009,
  title = {Statistical physics of social dynamics},
  author = {Castellano, Claudio and Fortunato, Santo and Loreto, Vittorio},
  journal = {Rev. Mod. Phys.},
  volume = {81},
  issue = {2},
  pages = {591--646},
  numpages = {0},
  year = {2009},
  month = {May},
  publisher = {American Physical Society},
  doi = {10.1103/RevModPhys.81.591},
  url = {https://link.aps.org/doi/10.1103/RevModPhys.81.591}
}

@article{DODDS_2005,
title = {A generalized model of social and biological contagion},
journal = {Journal of Theoretical Biology},
volume = {232},
number = {4},
pages = {587-604},
year = {2005},
issn = {0022-5193},
doi = {https://doi.org/10.1016/j.jtbi.2004.09.006},
url = {https://www.sciencedirect.com/science/article/pii/S0022519304004515},
author = {P.S. Dodds and D.J. Watts},
}

@article{BUMA_2020,
  title={Disturbances as opportunities: Learning from disturbance-response parallels in social and ecological systems to better adapt to climate change},
  author={Buma, Brian and Schultz, Courtney},
  journal={Journal of Applied Ecology},
  volume={57},
  number={6},
  pages={1113--1123},
  year={2020},
  publisher={Wiley Online Library}
}

@article{KEFI_2019,
  title={Advancing our understanding of ecological stability},
  author={K{\'e}fi, Sonia and Dom{\'\i}nguez-Garc{\'\i}a, Virginia and Donohue, Ian and Fontaine, Colin and Th{\'e}bault, Elisa and Dakos, Vasilis},
  journal={Ecology letters},
  volume={22},
  number={9},
  pages={1349--1356},
  year={2019},
  publisher={Wiley Online Library}
}

@ARTICLE{ZHI_2006,
  author={Zhi Li and Guanrong Chen},
  journal={IEEE Transactions on Circuits and Systems II: Express Briefs}, 
  title={Global synchronization and asymptotic stability of complex dynamical networks}, 
  year={2006},
  volume={53},
  number={1},
  pages={28-33},
  doi={10.1109/TCSII.2005.854315}}

@article{NEUBERT_1997,
  title={Alternatives to resilience for measuring the responses of ecological systems to perturbations},
  author={Neubert, Michael G and Caswell, Hal},
  journal={Ecology},
  volume={78},
  number={3},
  pages={653--665},
  year={1997},
  publisher={Wiley Online Library}
}

@article{PASCUAL_2005,
  title={Criticality and disturbance in spatial ecological systems},
  author={Pascual, Mercedes and Guichard, Fr{\'e}d{\'e}ric},
  journal={Trends in ecology \& evolution},
  volume={20},
  number={2},
  pages={88--95},
  year={2005},
  publisher={Elsevier}
}

@article{SAUNDERS_1975,
  title={Stability of complex ecosystems},
  author={Saunders, PT and Bazin, MJ},
  journal={Nature},
  volume={256},
  number={5513},
  pages={120--121},
  year={1975},
  publisher={Nature Publishing Group UK London}
}

@article{SCHOON_2012,
  title={Understanding disturbances and responses in social-ecological systems},
  author={Schoon, Michael L and Cox, Michael E},
  journal={Society \& Natural Resources},
  volume={25},
  number={2},
  pages={141--155},
  year={2012},
  publisher={Taylor \& Francis}
}

@book{STROGATZ_2018,
  author    = {Steven H. Strogatz},
  title     = {Nonlinear Dynamics and Chaos: With Applications to Physics, Biology, Chemistry, and Engineering},
  edition    = {2nd},
  year       = {2018},
  publisher  = {CRC Press},
  address    = {Boca Raton, FL},
  isbn       = {9780813349107}
}

@book{LYAPUNOV_ORIGINAL,
  title={The General Problem of the Stability of Motion},
  author={Lyapunov, Aleksandr Mikhailovich},
  translator={Fuller, A. T.},
  year={1992},
  publisher={Taylor \& Francis},
  address={London},
  note={Includes a biography by V. I. Smirnov and a bibliography of Lyapunov's work}
}

@article{SCHEFFER_2009,
  author={Scheffer, Marten and Bascompte, Jordi and Brock, William A. and Brovkin, Victor and Carpenter, Stephen R. and Dakos, Vasilis and Held, Hermann and van Nes, Egbert H. and Rietkerk, Max and Sugihara, George},
  title={Early-warning signals for critical transitions},
  journal={Nature},
  volume={461},
  number={7260},
  pages={53--59},
  year={2009},
  month={Sep},
  doi={10.1038/nature08227},
  url={https://doi.org/10.1038/nature08227},
  issn={1476-4687},
}

@article{MONTOYA_2006,
  author={Montoya, Jos{\'e} M. and Pimm, Stuart L. and Sol{\'e}, Ricard V.},
  title={Ecological networks and their fragility},
  journal={Nature},
  volume={442},
  number={7100},
  pages={259--264},
  year={2006},
  month={Jul},
  doi={10.1038/nature04927},
  pmid={16855581},
  issn={1476-4687},
}

@article{SUWEIS_2015,
  author  = {Samir Suweis and Jacopo Grilli and Jayanth R. Banavar and Stefano Allesina and Amos Maritan},
  title   = {Effect of localization on the stability of mutualistic ecological networks},
  journal = {Nature Communications},
  year    = {2015},
  volume  = {6},
  number  = {1},
  pages   = {10179},
  doi     = {10.1038/ncomms10179},
  url     = {https://doi.org/10.1038/ncomms10179},
  issn    = {2041-1723}
}

@book{D_J_WATTS_BOOK_2003,
  author = {D. J. Watts},
  title = {Six Degrees: The Science of a Connected Age},
  publisher = {W. W. Norton \&
Company, 2003.},
  year = {2003}
}

@article{
BARABASI_1999,
author = {Albert-László Barabási  and Réka Albert },
title = {Emergence of Scaling in Random Networks},
journal = {Science},
volume = {286},
number = {5439},
pages = {509-512},
year = {1999},
doi = {10.1126/science.286.5439.509},
URL = {https://www.science.org/doi/abs/10.1126/science.286.5439.509},
}

@book{BARABASI_2016,
  title = {Network Science},
  author = {Barabási, Albert-László},
  year = {2016},
  publisher = {Cambridge University Press}
}

@book{POINCARRE_1892,
  author    = {Poincar{\'{e}}, Henri},
  title     = {Les M{\'{e}}thodes Nouvelles de la M{\'{e}}canique C{\'{e}}leste},
  publisher = {Gauthier-Villars},
  address   = {Paris},
  year      = {1892},
  volume    = {1}
}

@book{ROUTH_1877,
  author    = {Routh, Edward John},
  title     = {A Treatise on the Stability of a Given State of Motion: Particularly Steady Motion},
  year      = {1877},
  publisher = {Macmillan and Company},
  address   = {London}
}

@article{HURWITZ_1895,
  author  = {Hurwitz, Adolf},
  title   = {On the Conditions under Which an Equation Has Only Roots with Negative Real Parts},
  journal = {Mathematische Annalen},
  volume  = {46},
  number  = {2},
  pages   = {273--284},
  year    = {1895},
  doi     = {10.1007/BF01446812}
}

@article{PERRON_1929,
  author  = {Perron, Oskar},
  title   = {Die Stabilit\"atsfrage bei Differentialgleichungen},
  journal = {Mathematische Zeitschrift},
  volume  = {32},
  pages   = {703--728},
  year    = {1929},
  doi     = {10.1007/BF01194662}
}

@article{GROBMAN_1959,
  author  = {Grobman, D. M.},
  title   = {Homeomorphism of Systems of Differential Equations},
  journal = {Doklady Akademii Nauk SSSR},
  volume  = {128},
  pages   = {880--881},
  year    = {1959}
}

@article{HARTMAN_1960,
  author  = {Hartman, Philip},
  title   = {A Lemma in the Theory of Structural Stability of Differential Equations},
  journal = {Proceedings of the American Mathematical Society},
  volume  = {11},
  number  = {4},
  pages   = {610--620},
  year    = {1960},
  doi     = {10.2307/2034715}
}

@article{KALMAN_1963,
  author  = {Kalman, Rudolf E.},
  title   = {Mathematical Description of Linear Dynamical Systems},
  journal = {Journal of the Society for Industrial and Applied Mathematics, Series A: Control},
  volume  = {1},
  number  = {2},
  pages   = {152--192},
  year    = {1963},
  doi     = {10.1137/0301010}
}

@article{GARDNER_ASHBY_1970,
  author  = {Gardner, Mark R. and Ashby, W. Ross},
  title   = {Connectance of Large Dynamic (Cybernetic) Systems: Critical Values for Stability},
  journal = {Nature},
  volume  = {228},
  number  = {5273},
  pages   = {784--784},
  year    = {1970},
  month   = {nov},
  doi     = {10.1038/228784a0},
  url     = {https://doi.org/10.1038/228784a0},
  issn    = {1476-4687}
}

@article{MAY_1972,
  author  = {May, Robert M.},
  title   = {Will a Large Complex System be Stable?},
  journal = {Nature},
  volume  = {238},
  number  = {5364},
  pages   = {413--414},
  year    = {1972},
  month   = {aug},
  doi     = {10.1038/238413a0},
  url     = {https://doi.org/10.1038/238413a0},
  issn    = {1476-4687}
}

@article{ALLESINA_2012,
  author  = {Allesina, Stefano and Tang, Si},
  title   = {Stability Criteria for Complex Ecosystems},
  journal = {Nature},
  volume  = {483},
  number  = {7388},
  pages   = {205--208},
  year    = {2012},
  month   = {mar},
  doi     = {10.1038/nature10832},
  url     = {https://doi.org/10.1038/nature10832},
  issn    = {1476-4687}
}

@article{TANG_PAWAR_ALLESINA_2014,
  author  = {Tang, Si and Pawar, Samraat and Allesina, Stefano},
  title   = {Correlation between Interaction Strengths Drives Stability in Large Ecological Networks},
  journal = {Ecology Letters},
  volume  = {17},
  number  = {9},
  pages   = {1094--1100},
  year    = {2014},
  month   = {sep},
  doi     = {10.1111/ele.12312},
  url     = {https://doi.org/10.1111/ele.12312},
  issn    = {1461-0248}
}

@article{SINHA_2005,
  title = {Evidence of universality for the May-Wigner stability theorem for random networks with local dynamics},
  author = {Sinha, Sitabhra and Sinha, Sudeshna},
  journal = {Phys. Rev. E},
  volume = {71},
  issue = {2},
  pages = {020902(R)},
  numpages = {4},
  year = {2005},
  month = {Feb},
  publisher = {American Physical Society},
  doi = {10.1103/PhysRevE.71.020902},
  url = {https://link.aps.org/doi/10.1103/PhysRevE.71.020902}
}

@article{KIRK_2015,
  author  = {Kirk, Paul and Rolando, Delphine M. Y. and MacLean, Adam L. and Stumpf, Michael P. H.},
  title   = {Conditional Random Matrix Ensembles and the Stability of Dynamical Systems},
  journal = {New Journal of Physics},
  volume  = {17},
  number  = {8},
  pages   = {083025},
  year    = {2015},
  month   = {aug},
  doi     = {10.1088/1367-2630/17/8/083025},
  url     = {https://doi.org/10.1088/1367-2630/17/8/083025},
  publisher = {IOP Publishing and Deutsche Physikalische Gesellschaft}
}

@article{TARNOWSKI_2020,
  title = {Universal transient behavior in large dynamical systems on networks},
  author = {Tarnowski, Wojciech and Neri, Izaak and Vivo, Pierpaolo},
  journal = {Phys. Rev. Res.},
  volume = {2},
  issue = {2},
  pages = {023333},
  numpages = {18},
  year = {2020},
  month = {Jun},
  publisher = {American Physical Society},
  doi = {10.1103/PhysRevResearch.2.023333},
  url = {https://link.aps.org/doi/10.1103/PhysRevResearch.2.023333}
}

@article{OSULLIVAN_2019,
  author  = {O'Sullivan, Jacob D. and Knell, Robert J. and Rossberg, Axel G.},
  title   = {Metacommunity-scale Biodiversity Regulation and the Self-Organised Emergence of Macroecological Patterns},
  journal = {Ecology Letters},
  volume  = {22},
  number  = {9},
  pages   = {1428--1438},
  year    = {2019},
  month   = {sep},
  doi     = {10.1111/ele.13294},
  url     = {https://doi.org/10.1111/ele.13294},
  issn    = {1461-0248}
}

@article{BARBIER_2021,
  title = {Fingerprints of High-Dimensional Coexistence in Complex Ecosystems},
  author = {Barbier, Matthieu and de Mazancourt, Claire and Loreau, Michel and Bunin, Guy},
  journal = {Phys. Rev. X},
  volume = {11},
  issue = {1},
  pages = {011009},
  numpages = {16},
  year = {2021},
  month = {Jan},
  publisher = {American Physical Society},
  doi = {10.1103/PhysRevX.11.011009},
  url = {https://link.aps.org/doi/10.1103/PhysRevX.11.011009}
}

@article{SILJAK_1973,
  author={Grujic, L. and Siljak, D.},
  journal={IEEE Transactions on Automatic Control}, 
  title={Asymptotic stability and instability of large-scale systems}, 
  year={1973},
  volume={18},
  number={6},
  pages={636-645},
  doi={10.1109/TAC.1973.1100422}}

@book{SILJAK_BOOK_1978,
  author    = {Siljak, Dragoslav D.},
  title     = {Large-Scale Dynamic Systems: Stability and Structure},
  publisher = {North-Holland},
  address   = {New York},
  year      = {1978},
  series    = {North-Holland Series in System Science and Engineering}
}

@ARTICLE{SILVA_2024,
  author={Silva, Guilherme Fróes and Donaire, Alejandro and Middleton, Richard and McFadyen, Aaron and Ford, Jason},
  journal={IEEE Transactions on Automatic Control}, 
  title={Scalable Input-to-State Stability of Nonlinear Interconnected Systems}, 
  year={2025},
  volume={70},
  number={3},
  pages={1824-1834},
  doi={10.1109/TAC.2024.3468069}}

@article{SONTANG_1989,
  title={Smooth stabilization implies coprime factorization},
  author={Sontag, Eduardo D.},
  journal={IEEE Transactions on Automatic Control},
  volume={34},
  number={4},
  pages={435--443},
  year={1989}
}

@article{MIRONCHENKO_2020,
author = {Mironchenko, Andrii and Prieur, Christophe},
title = {Input-to-State Stability of Infinite-Dimensional Systems: Recent Results and Open Questions},
journal = {SIAM Review},
volume = {62},
number = {3},
pages = {529-614},
year = {2020},
doi = {10.1137/19M1291248},
URL = {https://doi.org/10.1137/19M1291248},
}

@ARTICLE{TANG_2026,
  author={Tang, Michael and Krstić, Miroslav and Poveda, Jorge I.},
  journal={IEEE Transactions on Automatic Control}, 
  title={Fixed-Time Input-to-State Stability for Singularly Perturbed Systems via Composite Lyapunov Functions}, 
  year={2026},
  volume={71},
  number={4},
  pages={2644-2651},
  doi={10.1109/TAC.2025.3631559}}

@article{MEENA_2023,
  author  = {Meena, Chandrakala and Hens, Chittaranjan and Acharyya, Suman and Haber, Simcha and Boccaletti, Stefano and Barzel, Baruch},
  title   = {Emergent Stability in Complex Network Dynamics},
  journal = {Nature Physics},
  volume  = {19},
  number  = {7},
  pages   = {1033--1042},
  year    = {2023},
  month   = {jul},
  doi     = {10.1038/s41567-023-02020-8},
  url     = {https://doi.org/10.1038/s41567-023-02020-8},
  issn    = {1745-2481}
}

@article{BARZEL_2013,
  author  = {Barzel, Baruch and Barab{\'a}si, Albert-L{\'a}szl{\'o}},
  title   = {Universality in Network Dynamics},
  journal = {Nature Physics},
  volume  = {9},
  number  = {10},
  pages   = {673--681},
  year    = {2013},
  month   = {oct},
  doi     = {10.1038/nphys2741},
  url     = {https://doi.org/10.1038/nphys2741},
  issn    = {1745-2481}
}

@article{
HUFNAGEL_2004,
author = {L. Hufnagel  and D. Brockmann  and T. Geisel },
title = {Forecast and control of epidemics in a globalized world},
journal = {Proceedings of the National Academy of Sciences},
volume = {101},
number = {42},
pages = {15124-15129},
year = {2004},
doi = {10.1073/pnas.0308344101},
URL = {https://www.pnas.org/doi/abs/10.1073/pnas.0308344101},
}

@article{PASTOR_2001,
  title = {Epidemic Spreading in Scale-Free Networks},
  author = {Pastor-Satorras, Romualdo and Vespignani, Alessandro},
  journal = {Phys. Rev. Lett.},
  volume = {86},
  issue = {14},
  pages = {3200--3203},
  numpages = {0},
  year = {2001},
  month = {Apr},
  publisher = {American Physical Society},
  doi = {10.1103/PhysRevLett.86.3200},
  url = {https://link.aps.org/doi/10.1103/PhysRevLett.86.3200}
}

@book{VOIT_2000,
  author    = {Voit, Eberhard O.},
  title     = {Computational Analysis of Biochemical Systems},
  publisher = {Cambridge University Press},
  address   = {Cambridge},
  year      = {2000},
  isbn      = {9780521780223}
}

@article{MASLOV_2007,
  author  = {Maslov, Sergei and Ispolatov, Iaroslav},
  title   = {Propagation of Large Concentration Changes in Reversible Protein-Binding Networks},
  journal = {Proceedings of the National Academy of Sciences of the United States of America},
  volume  = {104},
  number  = {34},
  pages   = {13655--13660},
  year    = {2007},
  doi     = {10.1073/pnas.0704307104},
  url     = {https://doi.org/10.1073/pnas.0704307104}
}

@article{MASLOV_SERGEI_2007,
  author  = {Maslov, Sergei and Ispolatov, Iaroslav},
  title   = {Spreading Out of Perturbations in Reversible Reaction Networks},
  journal = {New Journal of Physics},
  volume  = {9},
  pages   = {273},
  year    = {2007},
  doi     = {10.1088/1367-2630/9/8/273},
  url     = {https://doi.org/10.1088/1367-2630/9/8/273}
}

@article{YAN_MASLOV_2008,
  author  = {Yan, Koon-Kiu and Walker, David and Maslov, Sergei},
  title   = {Fluctuations in Mass-Action Equilibrium of Protein Binding Networks},
  journal = {Physical Review Letters},
  volume  = {101},
  number  = {26},
  pages   = {268102},
  year    = {2008},
  doi     = {10.1103/PhysRevLett.101.268102},
  url     = {https://doi.org/10.1103/PhysRevLett.101.268102}
}

@book{GARDINER_2004,
  author    = {Gardiner, Crispin W.},
  title     = {Handbook of Stochastic Methods: For Physics, Chemistry and the Natural Sciences},
  publisher = {Springer},
  address   = {Berlin},
  edition   = {3rd},
  year      = {2004},
  isbn      = {9783540208822},
  doi       = {10.1007/978-3-662-05389-8}
}

@book{GANESH_2004,
  author    = {Hayes, John F. and Ganesh Babu, T. V. J.},
  title     = {Modeling and Analysis of Telecommunications Networks},
  publisher = {John Wiley \& Sons},
  address   = {Hoboken, New Jersey},
  year      = {2004},
  isbn      = {9780471471329}
}

@article{NOVOZHILOV_2006,
  author  = {Novozhilov, Artem S. and Karev, Georgy P. and Koonin, Eugene V.},
  title   = {Biological Applications of the Theory of Birth-and-Death Processes},
  journal = {Briefings in Bioinformatics},
  volume  = {7},
  number  = {1},
  pages   = {70--85},
  year    = {2006},
  doi     = {10.1093/bib/bbk006},
  url     = {https://doi.org/10.1093/bib/bbk006}
}

@book{ALON_BOOK_2006,
  author    = {Alon, Uri},
  title     = {An Introduction to Systems Biology: Design Principles of Biological Circuits},
  publisher = {Chapman \& Hall/CRC},
  address   = {Boca Raton, Florida},
  year      = {2006},
  isbn      = {9781584886426}
}

@article{KARLEBACH_2008,
  author  = {Karlebach, Guy and Shamir, Ron},
  title   = {Modelling and Analysis of Gene Regulatory Networks},
  journal = {Nature Reviews Molecular Cell Biology},
  volume  = {9},
  number  = {10},
  pages   = {770--780},
  year    = {2008},
  doi     = {10.1038/nrm2503},
  url     = {https://doi.org/10.1038/nrm2503}
}

@book{NEWMAN_BOOK_2018,
  author    = {Newman, Mark E. J.},
  title     = {Networks},
  edition   = {Second},
  publisher = {Oxford University Press},
  address   = {Oxford, United Kingdom},
  year      = {2018},
  isbn      = {9780198805090}
}

@article{PASTOR_2015,
  title = {Epidemic processes in complex networks},
  author = {Pastor-Satorras, Romualdo and Castellano, Claudio and Van Mieghem, Piet and Vespignani, Alessandro},
  journal = {Rev. Mod. Phys.},
  volume = {87},
  issue = {3},
  pages = {925--979},
  numpages = {55},
  year = {2015},
  month = {Aug},
  publisher = {American Physical Society},
  doi = {10.1103/RevModPhys.87.925},
  url = {https://link.aps.org/doi/10.1103/RevModPhys.87.925}
}

@article{
DODDS_2003,
author = {Peter Sheridan Dodds  and Roby Muhamad  and Duncan J. Watts },
title = {An Experimental Study of Search in Global Social Networks},
journal = {Science},
volume = {301},
number = {5634},
pages = {827-829},
year = {2003},
doi = {10.1126/science.1081058},
URL = {https://www.science.org/doi/abs/10.1126/science.1081058},
}

@incollection {BROCKMANN_2009,
    AUTHOR = {Brockmann, Dirk},
     TITLE = {Human mobility and spatial disease dynamics},
 BOOKTITLE = {Reviews of nonlinear dynamics and complexity. {V}ol. 2},
     PAGES = {1--24},
 PUBLISHER = {Wiley-VCH Verlag, Weinheim},
      YEAR = {2009},
      ISBN = {978-3-527-40850-4},
   MRCLASS = {92D30 (37N25 90B20)},
  MRNUMBER = {3307843},
}

@article{WODARZ_2002,
  title={The importance of lytic and nonlytic immune responses in viral infections},
  author={Wodarz, Dominik and Christensen, Jan Pravsgaard and Thomsen, Allan Randrup},
  journal={TRENDS in Immunology},
  volume={23},
  number={4},
  pages={194--200},
  year={2002},
  publisher={Elsevier}
}

@book{MURRAY_1989,
  author    = {Murray, James D.},
  title     = {Mathematical Biology},
  publisher = {Springer},
  year      = {1989}
}

@article{BARUCH_2012,
  title = {Stochastic analysis of complex reaction networks using binomial moment equations},
  author = {Barzel, Baruch and Biham, Ofer},
  journal = {Phys. Rev. E},
  volume = {86},
  issue = {3},
  pages = {031126},
  numpages = {13},
  year = {2012},
  month = {Sep},
  publisher = {American Physical Society},
  doi = {10.1103/PhysRevE.86.031126},
  url = {https://link.aps.org/doi/10.1103/PhysRevE.86.031126}
}

@article{HOLLING_1959,
  title={Some characteristics of simple types of predation and parasitism1},
  author={Holling, Crawford S},
  journal={The canadian entomologist},
  volume={91},
  number={7},
  pages={385--398},
  year={1959},
  publisher={Cambridge University Press}
}

@article{HOLLAND_2002,
author = {J. Nathaniel Holland  and Donald L. DeAngelis  and Judith L. Bronstein  and null null },
title = {Population Dynamics and Mutualism: Functional Responses of Benefits and Costs},
journal = {The American Naturalist},
volume = {159},
number = {3},
pages = {231-244},
year = {2002},
doi = {10.1086/338510},
URL = {https://www.journals.uchicago.edu/doi/abs/10.1086/338510}
}

@article{ZHANG_2024,
  author  = {Zhang, Xing and Li, Zhitao and Gao, Lixin},
  title   = {Stability Analysis of a SAIR Epidemic Model on Scale-Free Community Networks},
  journal = {Mathematical Biosciences and Engineering},
  volume  = {21},
  number  = {3},
  pages   = {4648--4668},
  year    = {2024},
  month   = {feb},
  doi     = {10.3934/mbe.2024204},
  url     = {https://doi.org/10.3934/mbe.2024204},
  pmid    = {38549343},
  issn    = {1547-1063}
}

@article{GAO_2016,
  author  = {Gao, Jianxi and Barzel, Baruch and Barab{\'a}si, Albert-L{\'a}szl{\'o}},
  title   = {Universal Resilience Patterns in Complex Networks},
  journal = {Nature},
  volume  = {530},
  number  = {7590},
  pages   = {307--312},
  year    = {2016},
  month   = {feb},
  doi     = {10.1038/nature16948},
  url     = {https://doi.org/10.1038/nature16948},
  issn    = {1476-4687}
}

@article{YAN_2017,
  title={Degree heterogeneity and stability of ecological networks},
  author={Yan, Gang and Martinez, Neo D and Liu, Yang-Yu},
  journal={Journal of the Royal Society Interface},
  volume={14},
  number={131},
  pages={20170189},
  year={2017}
}

@book{KHALIL_2002,
  author    = {Khalil, Hassan K.},
  title     = {Nonlinear Systems},
  edition   = {3},
  publisher = {Prentice Hall},
  address   = {Upper Saddle River, NJ},
  year      = {2002}
}

\end{document}